\documentclass[A4paper, oneside]{article}

\usepackage{amsthm}

\newtheorem{theorem}{Theorem}[section]
\newtheorem{corollary}[theorem]{Corollary}

\newtheorem{proposition}[theorem]{Proposition}
\newtheorem{conjecture}[theorem]{Conjecture}

\theoremstyle{definition}
\newtheorem{definition}[theorem]{Definition}

\theoremstyle{remark}
\newtheorem{remark}[theorem]{Remark}
\newtheorem{example}[theorem]{Example}

\usepackage{amsthm}           	    	              	
\usepackage{amsmath}            	    	            
\usepackage{amssymb}                	    	        
\usepackage{enumerate}									
	
\usepackage{csquotes} 		                          	

\usepackage[T1]{fontenc}                                

\usepackage[a4paper]{geometry}                          

\usepackage{tikz-cd}                                    
\usepackage[all]{xy}                                    
\usepackage{pgf,tikz}
\usepackage{mathrsfs}
\usetikzlibrary{arrows}
\usetikzlibrary{patterns}
\usepackage{hyperref}

\newcommand{\End}{\mathrm{End}}
\newcommand{\Filt}{\mathrm{Filt}}
\newcommand{\Hom}{\mathrm{Hom}}
\newcommand{\Ext}{\mathrm{Ext}}
\newcommand{\Fac}{\mathrm{Fac}}
\newcommand{\Sub}{\mathrm{Sub}}
\newcommand{\proj}{\mathrm{proj}}
\newcommand{\inj}{\mathrm{inj}}
\newcommand{\T}{\mathcal{T}}
\newcommand{\cS}{\mathcal{S}}
\newcommand{\F}{\mathcal{F}}
\newcommand{\C}{\mathcal{C}}
\newcommand{\D}{\mathfrak{D}}
\newcommand{\Ch}{\mathfrak{C}}
\newcommand{\PP}{\mathbf{P}}
\newcommand{\QQ}{\mathbf{Q}}

\newcommand{\coker}{\mbox{coker}}
\newcommand{\X}{\mathcal{X}}
\newcommand{\W}{\mathcal{W}}
\newcommand{\Y}{\mathcal{Y}}

\newcommand{\A}{\mathcal{A}}

\renewcommand{\mod}{\mathrm{mod}}

\newcommand{\add}{\mathrm{add}}

\renewcommand{\t}{\mathrm{t}}
\renewcommand{\k}{\mathbb{K}}

\renewcommand{\P}{\mathcal{P}}

\renewcommand{\ker}{\mbox{ker}}

\newcommand{\rep}[1]{%
  {%
    \tiny%
    \begin{matrix}%
      #1%
    \end{matrix}%
  }%
}

\begin{document}

\title{$\tau$-tilting theory -- an introduction}
\author{Hipolito Treffinger}
\date{\today}

\maketitle

\begin{abstract}
The notion of $\tau$-tilting theory was introduced by Adachi, Iyama and Reiten at the beginning of the last decade and quickly became one of the most active areas of research in the representation theory of finite dimensional algebras.
The aim of these notes is two-fold. 
On the one hand, we want to give a friendly introduction to $\tau$-tilting theory to anyone with a small background in representation theory. 
On the other, we want to fill the apparent gap for a survey on the subject by collecting in one place many of the most important results in $\tau$-tilting theory.
\end{abstract}

\section{Introduction}

The term \textit{$\tau$-tilting theory} was coined by Adachi, Iyama and Reiten in \cite{AIR} at the beginning of the 2010s.
In their paper, the authors created a fresh approach to the study of two classical branches of the representation theory of finite dimensional algebras, namely tilting theory and Auslander-Reiten theory. 
The combination of these two subjects is clearly reflected in the name of this novel theory, where the greek letter $\tau$ represents the Auslander-Reiten translation in the module category of an algebra while the reference to tilting theory is obvious.

Our primary aim with these notes\footnote{This is a revised and extended version of the Lecture notes written for the LMS Autumn School in Algebra 2020. For more information visit \url{https://www.icms.org.uk/events/event/?id=1073}.} is to give a friendly introduction to $\tau$-tilting theory from a representation theoretic perspective. 
Here you can find a compilation of many important results on the subject, giving a special emphasis in the close relation between $\tau$-tilting theory and torsion theories. 
Some background in representation theory and Auslander-Reiten theory is desirable but not necessary to follow this exposition. 
We note that given the immense amount of work that has been done in the last decade on $\tau$-tilting theory, this is not a complete survey on the topic.
For instance, we do not cover the rich connections that $\tau$-tilting theory has with other branches of mathematics, such as combinatorics or algebraic geometry.

This article can be divided into five different parts. 
In the first part, that consists only of Section~\ref{sec:history}, we give a brief historic account of the events leading to the rise of $\tau$-tilting theory in representation theory. 
We note that this section is not necessary for the understanding of the rest of the paper and can be skipped. 
In the second, which consists of Section~\ref{sec:definitions} and Section~\ref{sec:torsionclasses}, we give the basic definitions of the theory. 
We also describe some of the different forms that $\tau$-tilting can adopt, namely support $\tau$-tilting modules, $\tau$-tilting pairs, functorially finite torsion classes, and $2$-term silting complexes. 
The third part, corresponding to Sections~\ref{sec:tiltingtheorem}--\ref{sec:tau-tiltingfinite}, is dedicated to compiling general results on the subject, including the so-called $\tau$-tilting reduction and the characterisation of $\tau$-tilting finite algebras.
The $K$-theory of $\tau$-tilting theory is discussed in the fourth part of the paper, namely Sections~\ref{sec:vectors} to \ref{sec:g-vectors}.
We finish these notes in Section~\ref{sec:wallandchamberstructure}, where we show how we can associate to every algebra a geometric object known as its \textit{wall-and-chamber structure} and we explain how this invariant encodes much of its $\tau$-tilting theory.

We warn the reader that we do not include in these notes any proofs; these can be found in the references given by each result. 
Given the short time that has passed since the introduction of $\tau$-tilting theory, to our knowledge, there is not much material on the subject available other than the original research papers, with the exception of \cite{Iyama2014}. 
For background material in representation theory, we recommend the textbooks \cite{AsSS, ARS, Schifflerbook}.
For survey materials on more classical tilting theory, the reader is encouraged to see \cite{Assem1990, HandbookAHK}.

\subsection{Notation}
In these notes $A$ is always a basic finite dimensional algebra over a field $\k$ that we assume is algebraically closed. 
For us $\mod\, A$ is the category of finitely presented right $A$-modules and $\tau$ denotes the Auslander-Reiten translation in $\mod\, A$. 

Given any $A$-module $M$, we denote by $|M|$ the number of isomorphism classes of indecomposable direct summands of $M$.
Throughout this document we assume that $n$ is the number of isomorphism classes of simple $A$-modules. 
Note that in this case $|A|=n$, since $A = \bigoplus_{i=1}^n P(i)$, where $P(i)$ denotes the $i$-th indecomposable projective module.

Also, unless otherwise specified, every module is assumed to be basic, meaning that the indecomposable direct summands of $M$ are pairwise non-isomorphic.

\subsection*{Acknowledgements}
These notes were written as a support for a series of three lectures entitled \textit{$\tau$-tilting theory for finite-dimensional algebras} framed in the \textbf{LMS Autumn Algebra School 2020} organised by David Jordan, Nadia Mazza and Sibylle Schroll and funded by the London Mathematical Society. 
The author is grateful to the organisers for giving him the opportunity to present these lectures. 
He also thanks Bethany Marsh and Jenny August and the anonymous referee for their careful reading and their insightful comments in a previous versions of this document.
Finally, he acknowledges the financial support of the \textbf{Hausdorff Center of Mathematics}.
The author is also funded by the Deutsche Forschungsgemeinschaft (DFG, German Research Foundation) under Germany's Excellence Strategy Grant EXC-2047/1-390685813 and by the European Union’s Horizon 2020 research and innovation programme under the Marie Sklodowska-Curie grant agreement No 893654.


\section{Towards $\tau$-tilting theory}\label{sec:history}

It can be argued that the modern study of representation theory started with the parallel developments of almost split sequences by Auslander and Reiten \cite{Auslander1974, AR1, AR2} (see also \cite{ReitenLectureNotes}) and the theory of quiver representations by Gabriel \cite{Gabriel1972, Gabriel1973}.
Gabriel showed two very important results using quivers. 
One of these results says that the representation theory of every finite dimensional algebra over an algebraically closed field can be understood using quiver representations. 
The formal statement is the following.

\begin{theorem}\label{thm:Gabriel1}\cite{Gabriel1972}
Let $A$ be a finite dimensional algebra over an algebraically closed field $\k$. 
Then $A$ is Morita equivalent to the algebra $\k Q/I$, the path algebra of the quiver $Q$ bounded by an adimissible ideal of relations $I$.
Moreover the quiver $Q$ is uniquely determined by $A$. 
\end{theorem}

In the literature, people refer to the quiver $Q$ determined by the algebra as the \textit{ordinary quiver} or the \textit{Gabriel quiver} or simply the \textit{quiver} of the algebra.
In these notes we take the latter option.
The reason to give it such names is that one can associate to each finite dimensional algebra another quiver known as the \textit{Auslander-Reiten quiver} of the algebra, which encodes all the almost split sequences in $\mod~A$.
For more information about the Auslander-Reiten theory of algebras, the reader is encouraged to see the course on this topic by Raquel Coelho-Simões in this same series. 
See also \cite[IV.5]{AsSS}.

The second result of Gabriel we want to mention here is the classification of hereditary algebras of finite representation type by means of Dynkin diagrams as follows.

\begin{theorem}\label{thm:Gabriel2}\cite{Gabriel1972}
Let $A$ be a connected hereditary representation-finite finite dimensional algebra over an algebraically closed field $\k$. 
Then $A$ is Morita equivalent to $\k \vec{\Delta}$, where $\vec{\Delta}$ is a quiver whose underlying graph is a Dynkin diagram $\Delta$ of type $\mathbb{A}, \mathbb{D}$ or $\mathbb{E}$. 
Moreover there is a one-to-one correspondence between the indecomposable representations of $A$ and the positive roots of the root system associated to $\Delta$.
\end{theorem}

When this result appeared, it came as a great surprise since many fundamental properties of the path algebra of a quiver depend on the orientations of the arrows. 
For instance, if we start with two quivers $Q_1, Q_2$ that correspond to two different orientations of the same Dynkin diagram $\Delta$ then the path algebra $\k Q_1$ is in general not isomorphic to the path algebra $\k Q_2$, not even as vector spaces.
Hence, there was no reason to believe that the number of indecomposable representations should be the same. 

\begin{example}\label{ex:running}
For instance, take the algebras $A$ and $A'$ to be the path algebras of the quivers
$$Q_A=\xymatrix{1\ar[r] &  2\ar[r]& 3 } \qquad Q_{A'}=\xymatrix{1\ar[r] &  2& 3\ar[l] }$$ 
of type $\mathbb{A}_3$. 
A quick calculation shows that $\dim_\k A = 6$ while $\dim_\k A' = 5$.
The Auslander-Reiten quivers of $A$ and $A'$ can be found in Figure~\ref{fig:Ar-quiverA3} and Figure~\ref{fig:A'r-quiverA3}, respectively. 
Here the arrows correspond to the irreducible morphisms in the module category and the dashed lines correspond to the Auslander-Reiten translation.
In these figures we can see that the number of indecomposable representations of $A$ and $A'$ coincide.
\begin{figure}
    \centering
			\begin{tikzpicture}[line cap=round,line join=round ,x=2.0cm,y=1.8cm]
				\clip(-2.2,-0.1) rectangle (4.1,2.5);
					\draw [->] (-0.8,0.2) -- (-0.2,0.8);
					\draw [->] (1.2,0.2) -- (1.8,0.8);
					\draw [->] (0.2,1.2) -- (0.8,1.8);
					\draw [<-, dashed] (-0.8,0.0) -- (0.8,0.0);
					\draw [<-, dashed] (1.2,0.0) -- (2.8,0.0);
					\draw [<-, dashed] (0.2,1.0) -- (1.8,1.0);
					\draw [->] (0.2,0.8) -- (0.8,0.2);
					\draw [->] (2.2,0.8) -- (2.8,0.2);
					\draw [->] (1.2,1.8) -- (1.8,1.2);
				
				\begin{scriptsize}
					\draw[color=black] (-1,0) node {$\rep{3}$};
					\draw[color=black] (1,0) node {$\rep{2}$};
					\draw[color=black] (3,0) node {$\rep{1}$};
					\draw[color=black] (0,1) node {$\rep{2\\3}$};
					\draw[color=black] (2,1) node {$\rep{1\\2}$};
					\draw[color=black] (1,2) node {$\rep{1\\2\\3}$};
				\end{scriptsize}
			\end{tikzpicture}
\caption{The Auslander-Reiten quiver of $A$}
    \label{fig:Ar-quiverA3}
\end{figure}
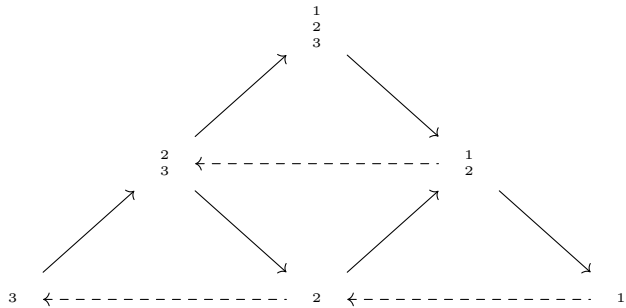

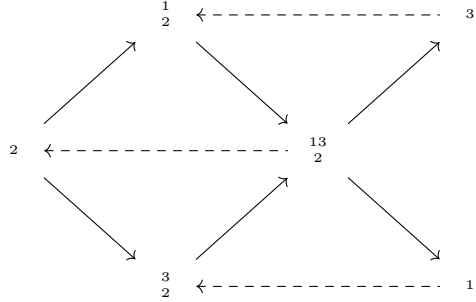
\begin{figure}
    \centering
			\begin{tikzpicture}[line cap=round,line join=round ,x=2.0cm,y=1.8cm]
				\clip(-1.2,-0.1) rectangle (4.1,2.5);
					\draw [->] (1.2,0.2) -- (1.8,0.8);
					\draw [->] (0.2,1.2) -- (0.8,1.8);
					\draw [->] (2.2,1.2) -- (2.8,1.8);
					\draw [<-, dashed] (1.2,0.0) -- (2.8,0.0);
					\draw [<-, dashed] (1.2,2.0) -- (2.8,2.0);
					\draw [<-, dashed] (0.2,1.0) -- (1.8,1.0);
					\draw [->] (0.2,0.8) -- (0.8,0.2);
					\draw [->] (2.2,0.8) -- (2.8,0.2);
					\draw [->] (1.2,1.8) -- (1.8,1.2);
				
				\begin{scriptsize}
					\draw[color=black] (1,0) node {$\rep{3\\2}$};
					\draw[color=black] (3,0) node {$\rep{1}$};
					\draw[color=black] (0,1) node {$\rep{2}$};
					\draw[color=black] (2,1) node {$\rep{13\\2}$};
					\draw[color=black] (1,2) node {$\rep{1\\2}$};
					\draw[color=black] (3,2) node {$\rep{3}$};
				\end{scriptsize}
			\end{tikzpicture}
\caption{The Auslander-Reiten quiver of $A'$}
    \label{fig:A'r-quiverA3}
\end{figure}

\end{example}

As a consequence, explaining this phenomenon became of significant interest. 
The first explanation was given by Bernstein, Gelfand and Ponomarev in \cite{BGP} by constructing the so-called reflection functors. 

Let $Q$ a quiver of type $\Delta$ and denote by $Q_0$ the set of vertices of $Q$.
Since every Dynkin diagram is a tree, there is at least one vertex $x\in Q_0$ which is a sink, i.e. a vertex such that all the arrows incident to that vertex are incoming arrows. 
Now, we construct a quiver $Q_{A'}$ which is identical to $Q_A$, except for the fact that now the vertex $x$ is a source, which means that every arrow incident to $x$ is an outgoing arrow. 
One says that $Q_A$ and $Q_{A'}$ are reflections of each other at $x$.
In Example~\ref{ex:running}, $Q'$ is the reflection of $Q$ at the vertex $3$.
Then, Bernstein, Gelfand and Ponomarev showed the existence of functors, that they called reflection functors, between $\mod\; \k Q$ and $\mod\; \k Q'$ that induce a one-to-one correspondence between their indecomposable objects. 

Some years after that, Auslander, Platzeck and Reiten \cite{APR} realised that these functors were induced by a very specific object in $\mod\;\k Q$.
To be more precise, note that the simple module $S(x)$ associated to the vertex $x \in Q_0$ is projective and it is not injective.
This implies that the inverse Auslander-Reiten translation $\tau^{-1} S(x)$ of $S(x)$ is a non-zero indecomposable object of $\mod\; \k Q$.
Then they showed that the reflection functors described by Berstein, Gelfand and Ponomarev were equivalent to $\Hom_A (T, -) : \mod\, A \to \mod\, A'$, where $T$ is the module
\begin{equation}
T = \tau^{-1} S(x) \oplus \bigoplus_{x\neq y \in Q_0} P(y).
\label{eq:APR}
\end{equation}
Thus, $T$ is the direct sum $\tau^{-1} S(x)$ and the direct sum of all of the indecomposable projectives except $S(x)$. 
Moreover, they showed that $\k Q'$ is isomorphic to $\End_{\k Q} (T)^{op}$.
In particular, this approach allowed them to show the existence of reflection-like functors between the module category of any Artin algebra $A$ having a simple projective module and $\End_{A} (T)$, even when $A$ is not hereditary or even when $A$ is not the quotient of the path algebra of a quiver. 
Going once again to our running example, the module described by Auslander, Platzeck and Reiten in $\mod\, A$ is $T = \rep{2} \oplus \rep{1\\2\\3} \oplus \rep{2\\3}$

Some years later, Brenner and Butler went further and studied in \cite{BB} this phenomenon axiomatically. 
In this paper they introduce the notion of \textit{tilting modules} as follows. 

\begin{definition}\cite{BB}\label{def:tilting}
Let $A$ be an algebra and $T$ be an $A$-module. 
We say that $T$ is a tilting module if the following holds:
\begin{enumerate}
	\item[(i)] $\text{pd}_A T \leq 1$, the projective dimension of $T$ is at most $1$.
	\item[(ii)] $T$ is rigid, that is $\Ext^1_A(T,T)=0$.
	\item[(iii)] There exists a short exact sequence of the form
	$$0 \to A \to T' \to T'' \to 0$$
	where $T', T''$ are direct summands of direct sums of $T$.
\end{enumerate}
\end{definition}

In this paper they show that any tilting $A$-module $T$ acts as a sort of translator between the representation theory of $A$ and $B:=\End_A(T)^{op}$, the opposite of the endomorphism algebra of $T$.

\begin{figure}
    \centering
			\begin{tikzpicture}[line cap=round,line join=round ,x=2.0cm,y=1.8cm]
					\draw[-,rounded corners=0.1cm, thick, blue, fill=blue!20] 
					(0, 1.2)--(-0.2,1)--(0.8,0)--(1, -0.2)-- (3.0, -0.2) -- (3.2, 0) -- 
					(3.0, 0.2)-- (1, 2.3) -- (0, 1.2);
					\draw[-, rounded corners=0.1cm, thick, green, fill=green!20]
					(-1.1, -0.1)--(-1.2,0)--(-1, 0.2)--(-0.8, 0)--(-1, -0.2)--(-1.2, 0)
					--(-1.1, 0.1);

				\clip(-2.2,-0.3) rectangle (4.1,2.5);
					\draw [->] (-0.8,0.2) -- (-0.2,0.8);
					\draw [->] (1.2,0.2) -- (1.8,0.8);
					\draw [->] (0.2,1.2) -- (0.8,1.8);
					\draw [<-, dashed] (-0.8,0.0) -- (0.8,0.0);
					\draw [<-, dashed] (1.2,0.0) -- (2.8,0.0);
					\draw [<-, dashed] (0.2,1.0) -- (1.8,1.0);
					\draw [->] (0.2,0.8) -- (0.8,0.2);
					\draw [->] (2.2,0.8) -- (2.8,0.2);
					\draw [->] (1.2,1.8) -- (1.8,1.2);
				\begin{scriptsize}
					\draw[color=black] (-1,0) node {$\rep{3}$};
					\draw[color=black] (1,0) node {$\rep{2}$};
					\draw[color=black] (3,0) node {$\rep{1}$};
					\draw[color=black] (0,1) node {$\rep{2\\3}$};
					\draw[color=black] (2,1) node {$\rep{1\\2}$};
					\draw[color=black] (1,2) node {$\rep{1\\2\\3}$};
					\draw[color=blue] (2.05, 1.5) node {\huge{$\mathcal{T}$}};
					\draw[color=green] (-1, 0.4) node {\huge{$\mathcal{F}$}};
				\end{scriptsize}
			\end{tikzpicture}
			
			    \centering
			\begin{tikzpicture}[line cap=round,line join=round ,x=2.0cm,y=1.8cm]
				\clip(-1.4,-0.3) rectangle (4.0,2.8);
					\draw[-,rounded corners=0.1cm, thick, blue, fill=blue!20] 
					(0, 1.2)--(-0.2,1)--(0.8,0)--(1, -0.2)-- (3.0, -0.2) -- (3.2, 0) -- 
					(3.0, 0.2)-- (1, 2.3) -- (0, 1.2);
					\draw[-, rounded corners=0.1cm, thick, green, fill=green!20]
					(2.9, 1.9)--(2.8,2)--(3, 2.2)--(3.2, 2)--(3, 1.8)--(2.8, 2)
					--(2.9, 2.1);
					\draw [->] (1.2,0.2) -- (1.8,0.8);
					\draw [->] (0.2,1.2) -- (0.8,1.8);
					\draw [->] (2.2,1.2) -- (2.8,1.8);
					\draw [<-, dashed] (1.2,0.0) -- (2.8,0.0);
					\draw [<-, dashed] (1.2,2.0) -- (2.8,2.0);
					\draw [<-, dashed] (0.2,1.0) -- (1.8,1.0);
					\draw [->] (0.2,0.8) -- (0.8,0.2);
					\draw [->] (2.2,0.8) -- (2.8,0.2);
					\draw [->] (1.2,1.8) -- (1.8,1.2);

				\begin{scriptsize}
					\draw[color=black] (1,0) node {$\rep{3\\2}$};
					\draw[color=black] (3,0) node {$\rep{1}$};
					\draw[color=black] (0,1) node {$\rep{2}$};
					\draw[color=black] (2,1) node {$\rep{13\\2}$};
					\draw[color=black] (1,2) node {$\rep{1\\2}$};
					\draw[color=black] (3,2) node {$\rep{3}$};
					\draw[color=blue] (0, 1.5) node {\huge{$\mathcal{Y}$}};
					\draw[color=green] (3, 1.5) node {\huge{$\mathcal{X}$}};
				\end{scriptsize}
			\end{tikzpicture}
	\caption{Torsion pairs and APR-tilting}
    \label{fig:APRtilting}
\end{figure}
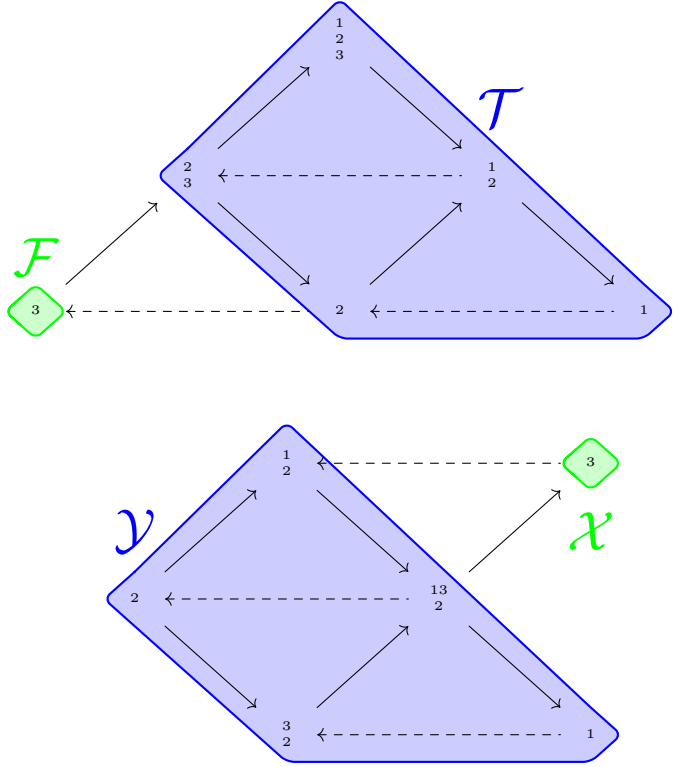

The first thing that they have shown is that a tilting $A$-module $T$ is also a tilting $B$-module.
Moreover they showed that $T$ induced a torsion pair $(\T, \F)$ in $\mod\, A$ and a torsion pair $(\X, \Y)$ in $\mod\, B$ such that the functors $\Hom_A(T,-): \mod\, A \to \mod\, B$ and $\Ext^1_A(T,-): \mod\, A \to \mod\, B$ induce equivalences of categories between $\T$ and $\Y$ and between $\F$ and $\X$, respectively.
This result of Brenner and Butler can be seen applied to our running example in Figure~\ref{fig:APRtilting}.
For the precise definition of torsion pair, see Definition~\ref{def:torsionpair}.
Also, a more detailed treatment of the tilting theorem will be given in Section~\ref{sec:tiltingtheorem}.

Since the module introduced by Auslander, Platzeck and Reiten was their motivating example, one can expect that it satisfies Definition~\ref{def:tilting} (i)-(iii) and indeed this is the case. 
In fact, nowadays this module is known as the \textit{APR-tilting module}.
But, as the reader is already guessing, there are many more examples of tilting modules.
Take the module $T = \rep{3} \oplus \rep{1\\2\\3} \oplus \rep{1}$.
One can verify that $T$ is indeed a tilting module. 
Firstly, the projective dimension of $T$ is less or equal to one since $A$ is hereditary. 
Secondly, one can check that $T$ does not admit self extensions. 
Finally, the short exact sequence 
$$0 \to \rep{3} \oplus \rep{2\\3} \oplus \rep{1\\2\\3} \to \rep{3} \oplus \rep{1\\2\\3} \oplus \rep{1\\2\\3} \to \rep{1} \to 0$$
is such that $\rep{3} \oplus \rep{1\\2\\3} \oplus \rep{1\\2\\3}$ and $\rep{1}$ are direct summands of direct sums of copies of $T$.

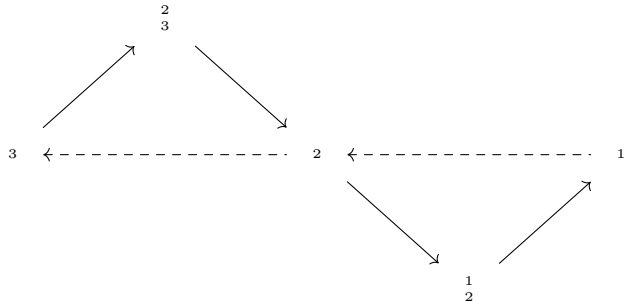
\begin{figure}
    \centering
			\begin{tikzpicture}[line cap=round,line join=round ,x=2.0cm,y=1.8cm]
				\clip(-0.4,-0.3) rectangle (4.3,2.5);
					\draw [->] (3.2,0.2) -- (3.8,0.8);
					\draw [->] (0.2,1.2) -- (0.8,1.8);
					\draw [<-, dashed] (0.2,1.0) -- (1.8,1.0);
					\draw [<-, dashed] (2.2,1.0) -- (3.8,1.0);
					\draw [->] (2.2,0.8) -- (2.8,0.2);
					\draw [->] (1.2,1.8) -- (1.8,1.2);
				
				\begin{scriptsize}
					\draw[color=black] (3,0) node {$\rep{1\\2}$};
					\draw[color=black] (4,1) node {$\rep{1}$};
					\draw[color=black] (0,1) node {$\rep{3}$};
					\draw[color=black] (2,1) node {$\rep{2}$};
					\draw[color=black] (1,2) node {$\rep{2\\3}$};
				\end{scriptsize}
			\end{tikzpicture}
\caption{The Auslander-Reiten quiver of $B$}
    \label{fig:Ar-quiverA3incl}
\end{figure}

Now, the algebra $B = \End_A(T)$ is isomorphic to the path algebra of the quiver  
$$\xymatrix{1\ar[r] &  2\ar[r]& 3}$$ 
modulo the ideal generated by the composition of the two arrows. 
The Auslander-Reiten quiver of $B$ can be seen in Figure~\ref{fig:Ar-quiverA3incl}.

\begin{figure}
    \centering
			\begin{tikzpicture}[line cap=round,line join=round ,x=2.0cm,y=1.8cm]
				\clip(-2.2,-0.3) rectangle (4.1,2.5);
					\draw[-,rounded corners=0.1cm, thick, blue, fill=blue!20] 
					(3.0, 0.2)--(1, 2.3)--(0.8,2)--(1.8,1)--(2, 0.8)-- 
					(3.0, -0.2) -- (3.2, 0) -- 
					(3.0, 0.2)-- (1, 2.3)--(0.8,2)--(1.8,1) ;
					\draw[-, rounded corners=0.1cm, thick, blue, fill=blue!20]
					(-1.1, -0.1)--(-1.2,0)--(-1, 0.2)--(-0.8, 0)--(-1, -0.2)--(-1.2, 0)
					--(-1.1, 0.1);
					\draw[-, rounded corners=0.1cm, thick, green, fill=green!20]
					(0.9, -0.1)--(0.8,0)--(1, 0.2)--(1.2, 0)--(1, -0.2)--(0.8, 0)
					--(0.9, 0.1);
					\draw [->] (-0.8,0.2) -- (-0.2,0.8);
					\draw [->] (1.2,0.2) -- (1.8,0.8);
					\draw [->] (0.2,1.2) -- (0.8,1.8);
					\draw [<-, dashed] (-0.8,0.0) -- (0.8,0.0);
					\draw [<-, dashed] (1.2,0.0) -- (2.8,0.0);
					\draw [<-, dashed] (0.2,1.0) -- (1.8,1.0);
					\draw [->] (0.2,0.8) -- (0.8,0.2);
					\draw [->] (2.2,0.8) -- (2.8,0.2);
					\draw [->] (1.2,1.8) -- (1.8,1.2);
				
				\begin{scriptsize}
					\draw[color=black] (-1,0) node {$\rep{3}$};
					\draw[color=black] (1,0) node {$\rep{2}$};
					\draw[color=black] (3,0) node {$\rep{1}$};
					\draw[color=black] (0,1) node {$\rep{2\\3}$};
					\draw[color=black] (2,1) node {$\rep{1\\2}$};
					\draw[color=black] (1,2) node {$\rep{1\\2\\3}$};
					\draw[color=blue] (2.2, 1.5) node {\huge{$\mathcal{T}$}};
					\draw[color=green] (1, 0.5) node {\huge{$\mathcal{F}$}};

				\end{scriptsize}
			\end{tikzpicture}

    \centering
			\begin{tikzpicture}[line cap=round,line join=round ,x=2.0cm,y=1.8cm]
				\clip(-0.4,-0.3) rectangle (4.3,3.5);
					\draw[-,rounded corners=0.1cm, thick, blue, fill=blue!20] 
					(1.2, 2)--(1, 2.2)--(0.8,2)--(1.8,1)--(2, 0.8)-- 
					(3.0, -0.2) -- (3.2, 0) -- 
					(4.0, 0.8)-- (4.2, 1) -- (4.0, 1.2) -- (2, 1.2)-- (1.2, 2);
					\draw[-, rounded corners=0.1cm, thick, green, fill=green!20]
					(-0.1, 0.9)--(-0.2,1)--(0, 1.2)--(0.2, 1)--(0, 0.8)--(-0.2, 1)
					--(-0.1, 1.1);
					\draw [->] (3.2,0.2) -- (3.8,0.8);
					\draw [->] (0.2,1.2) -- (0.8,1.8);
					\draw [<-, dashed] (0.2,1.0) -- (1.8,1.0);
					\draw [<-, dashed] (2.2,1.0) -- (3.8,1.0);
					\draw [->] (2.2,0.8) -- (2.8,0.2);
					\draw [->] (1.2,1.8) -- (1.8,1.2);
				
				\begin{scriptsize}
					\draw[color=black] (3,0) node {$\rep{1\\2}$};
					\draw[color=black] (4,1) node {$\rep{1}$};
					\draw[color=black] (0,1) node {$\rep{3}$};
					\draw[color=black] (2,1) node {$\rep{2}$};
					\draw[color=black] (1,2) node {$\rep{2\\3}$};
					\draw[color=blue] (2.2, 1.5) node {\huge{$\mathcal{Y}$}};
					\draw[color=green] (0, 1.5) node {\huge{$\mathcal{X}$}};
				\end{scriptsize}
			\end{tikzpicture}
\caption{The tilting theorem applied to a tilted algebra}
    \label{fig:torsionAB}
\end{figure}
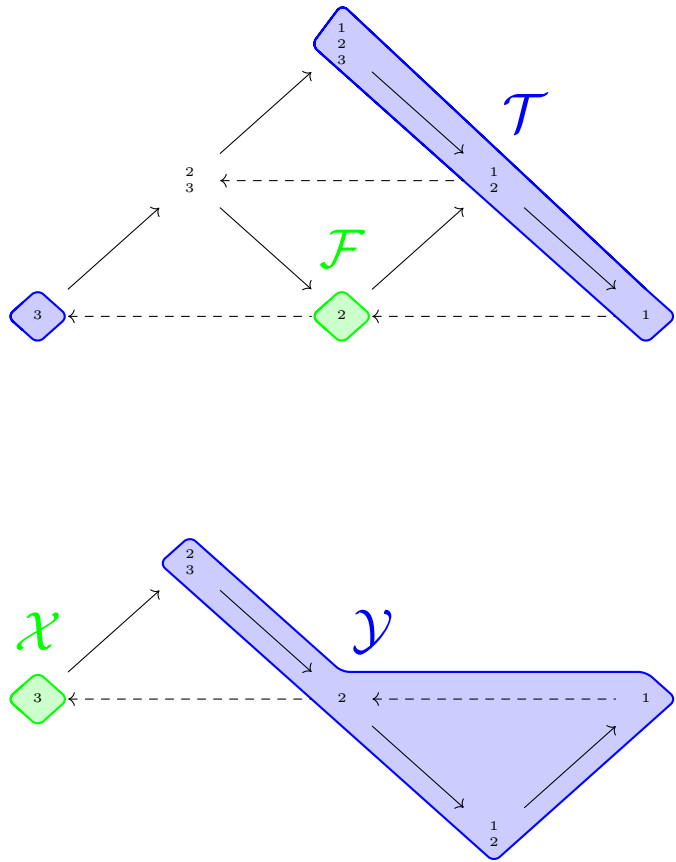

As we can see in this example, when we take an arbitrary tilting module $A$ the numbers of indecomposable representations in $\mod\, A$ and in $\mod\, B$ are not the same.
However, this is not a contradiction of the results of Brenner and Butler since their result only says what happens inside the torsion pairs induced by $T$ in $\mod\, A$ and $\mod\, B$. 
In this particular case, as we can see in Figure~\ref{fig:torsionAB}, the indecomposable object $\rep{2\\3}$ does not belong to either of the two subcategories $\T$ and $\F$ induced by the tilting module $T$.

Although the categories $\mod\, A$ and $\mod\,\End_A(T)$ are not in general equivalent, it was shown by Happel \cite{Happel1988} (see also Rickard's generalisation \cite{Rickard1989}) that these two algebras are \textit{derived equivalent}.
Without going into the details, starting from the module category of an algebra, one can construct a triangulated category known as the \textit{derived category} of the algebra that encodes a wealth of homological information of the algebra. 
Then, the results of Happel and Rickard state that the original algebra and the endomorphism algebra of the tilting module have the same derived category, which implies that they share many homological properties that we will not discuss here.
For more information about derived categories and their relationship with cluster algebras, we refer the reader to the notes of Matthew Pressland in this series. 

Let us include a small parenthesis here that will be important later. 
The algebra $B$ is the smallest non-hereditary example of a so-called \textit{tilted algebra}. 
Tilted algebras were introduced by Happel and Ringel in \cite{HR} (see also Bongartz \cite{BongartzComp}) as the endomorphism algebras of tilting modules over hereditary algebras. 
The main idea behind their introduction was to use all the information available on hereditary algebras to understand a new class of algebras which had not been studied systematically until that moment. 

The study of tilted algebras has sparked a great deal of research which it would be impossible to describe completely here. 
However, we need to mention two famous developments. 

Firstly, note that the tilting theorem of Brenner and Butler does not impose any restriction on the algebra $A$. 
So, if we are able to understand some of the representation theory of tilted algebras using the knowledge we have on hereditary algebras, we can repeat the process and understand the representation theory of a new family of algebras using the knowledge we have on tilted algebras via the tilting theorem. 
These algebras are known as \textit{iterated tilted algebras}. 
In \cite{AsSkoGentle}, Assem and Skowro\'nski  classified all the iterated tilted algebras of Dynkin type $\tilde{\mathbb{A}}$ in terms of their ordinary quiver and relations, which led them to the definition of the so-called \textit{gentle algebras}. 
Today, gentle algebras constitute a highly active area of research, deepening our understanding of representation theory of finite dimensional algebras and connecting this topic with various other branches of mathematics such as group theory and algebraic and differential geometry.

The second is the characterisation of tilted algebras found independently by Liu \cite{LiuCrit} and Skowro\'nski \cite{Skocrit} using the Auslander-Reiten quiver an algebra.
They have shown that an algebra is tilted if and only if there is a structure with specific homological and combinatorial properties in their Auslander-Reiten quiver. 
Inspired by this characterisation of tilted algebras many families of algebras have been defined and determined by means of their Auslander-Reiten quivers. 

\medskip
Some years later, at the beginning of the twenty-first century, Fomin and Zelevinsky \cite{FZ1, FZ2, BFZ3, FZ4} were studying the properties of the canonical bases arising in Lie theory and this study led to the introduction of \textit{cluster algebras}. 

These algebras are generated by a set of so-called \textit{cluster variables} that are produced inductively from an \textit{initial seed} via a process called \textit{mutation} that produces new seeds.
Even though the process of mutation is iterated an arbitrary (finite) number of times, for some initial seeds there are only finitely many cluster variables that can be constructed.
In this case we say that a cluster algebra is of \textit{finite type}.
Moreover, for some of these algebras, known as \textit{skew-symmetrizable} cluster algebras, their combinatorial construction can be expressed using quivers. 
One surprising result shown by Fomin and Zelevinsky in the first of the series of papers where they introduced cluster algebras is the following classification. 

\begin{theorem}\cite{FZ2}
Let $(Q, \{\underline{x}\})$ be the initial seed of a cluster algebra $\mathcal{A}$. 
Then $\mathcal{A}$ is of finite type if and only if $Q$ is mutation equivalent to a quiver whose underlying graph is a Dynkin diagram.
\end{theorem}

The resemblance of this result with Theorem~\ref{thm:Gabriel2} is striking and points towards a deep relationship between cluster theory and representations of finite dimensional algebras. 

It is very important to remark that for any seed $(Q', \{\underline{x'}\})$ the set of cluster variables $\{\underline{x}\}$ always has the same number of elements, lets call this number $n$.
Then, all the seeds of a cluster algebra can be arranged into a $n$-regular graph where there is an edge between two seeds if one can be obtained from the other performing a single mutation. 

As it turns out, similar phenomena have been described in tilting theory. 
For instance it was shown by Skowro\'nski in \cite{Skolemma} that every basic tilting module has exactly $n$ indecomposable direct summands. 
Also, Happel and Unger have shown in \cite{HU1} that every basic partial tilting module having $n-1$ indecomposable direct summands can always be completed into a tilting module and that there are at most two ways in which this can be done.

Hence, one would like to categorify all the cluster phenomena using tilting theory, where the cluster variables are represented by indecomposable partial tilting modules and tilting modules correspond to seeds.
However, tilting theory falls short in describing the cluster phenomena for at least two reasons. 
The first is that there are some examples of almost complete tilting modules that can be completed into a tilting module in exactly one way, which means that we can not reproduce the process of mutation at some indecomposable direct summand of this module. 

The second reason, and maybe the most obvious, is that there are fewer indecomposable partial tilting modules than cluster variables.
For instance, a hereditary path algebra of type $\mathbb{A}_n$ has exactly $\frac{n(n-1)}{2}$ indecomposable partial tilting modules, while the number of cluster variables in a cluster algebra of type $\mathbb{A}_n$ is $\frac{n(n+1)}{2}$, i.e. there are exactly $n$ more cluster variables than indecomposable partial tilting modules. 

Then if one wants to categorify cluster algebras using tilting theory, it is necessary to extend the latter in some way. 
That is exactly what Buan, Marsh, Reineke, Reiten and Todorov did in \cite{BMRRT}. 
In this seminal paper, instead of working with the module category of the algebra,  they constructed a slightly larger triangulated category that they called the \textit{cluster category} where everything works perfectly by the definition of the so-called \textit{(partial) cluster-tilting objects}.

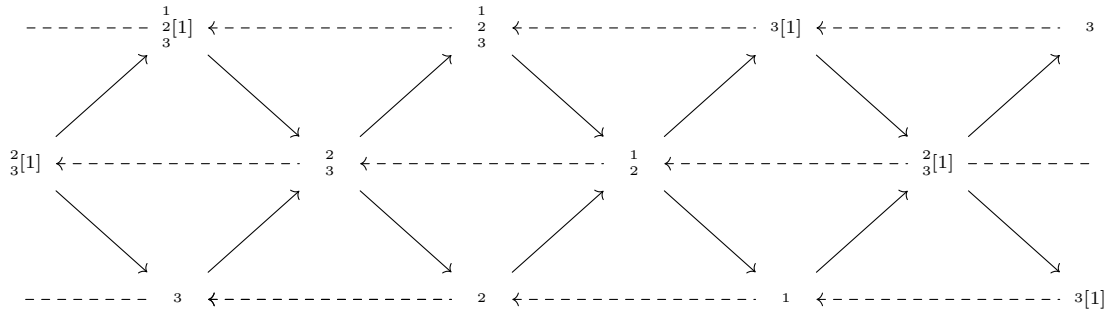
\begin{figure}
    \centering
			\begin{tikzpicture}[line cap=round,line join=round ,x=2.0cm,y=1.8cm]
				\clip(-2.2,-0.1) rectangle (5.1,2.5);
					\draw [->] (-0.8,0.2) -- (-0.2,0.8);
					\draw [->] (4.2,0.8) -- (4.8,0.2);
					\draw [->] (4.2,1.2) -- (4.8,1.8);
					\draw [->] (1.2,0.2) -- (1.8,0.8);
					\draw [->] (3.2,0.2) -- (3.8,0.8);
					\draw [->] (-1.8,1.2) -- (-1.2,1.8);
					\draw [->] (0.2,1.2) -- (0.8,1.8);
					\draw [->] (2.2,1.2) -- (2.8,1.8);
					\draw [<-, dashed] (-0.8,2.0) -- (0.8,2.0);
					\draw [<-, dashed] (1.2,2.0) -- (2.8,2.0);
					\draw [<-, dashed] (3.2,2.0) -- (4.8,2.0);
					\draw [dashed] (4.2,1.0) -- (5.0,1.0);
					\draw [dashed] (-1.2,2.0) -- (-2.0,2.0);
					\draw [<-, dashed] (-0.8,0.0) -- (0.8,0.0);
					\draw [<-, dashed] (-0.8,0.0) -- (0.8,0.0);
					\draw [<-, dashed] (1.2,0.0) -- (2.8,0.0);
					\draw [<-, dashed] (-1.8,1.0) -- (-0.2,1.0);
					\draw [<-, dashed] (0.2,1.0) -- (1.8,1.0);
					\draw [<-, dashed] (2.2,1.0) -- (3.8,1.0);
					\draw [dashed] (-2.0,0.0) -- (-1.2,0.0);
					\draw [<-, dashed] (3.2,0.0) -- (4.8,0.0);
					\draw [->] (0.2,0.8) -- (0.8,0.2);
					\draw [->] (2.2,0.8) -- (2.8,0.2);
					\draw [->] (-1.8,0.8) -- (-1.2,0.2);
					\draw [->] (-0.8,1.8) -- (-0.2,1.2);
					\draw [->] (1.2,1.8) -- (1.8,1.2);
					\draw [->] (3.2,1.8) -- (3.8,1.2);
				
				\begin{scriptsize}
					\draw[color=black] (-1,0) node {$\rep{3}$};
					\draw[color=black] (5,0) node {$\rep{3}[1]$};
					\draw[color=black] (5,2) node {$\rep{3}$};
					\draw[color=black] (1,0) node {$\rep{2}$};
					\draw[color=black] (3,0) node {$\rep{1}$};
					\draw[color=black] (-2,1) node {$\rep{2\\3}[1]$};
					\draw[color=black] (0,1) node {$\rep{2\\3}$};
					\draw[color=black] (2,1) node {$\rep{1\\2}$};
					\draw[color=black] (4,1) node {$\rep{2\\3}[1]$};
					\draw[color=black] (-1,2) node {$\rep{1\\2\\3}[1]$};
					\draw[color=black] (1,2) node {$\rep{1\\2\\3}$};
					\draw[color=black] (3,2) node {$\rep{3}[1]$};
				\end{scriptsize}
			\end{tikzpicture}
\caption{The Auslander-Reiten quiver of $\mathcal{C}_A$}
    \label{fig:Clustercat}
\end{figure}

See in Figure~\ref{fig:Clustercat} the Auslander-Reiten quiver of the cluster category associated to the algebra $A$ of Example~\ref{ex:running}. 
The points that are tagged with the same object in the Auslander-Reiten quiver of $\mathcal{C}_A$ should be identified. 
In particular, we see that the Auslander-Reiten quiver of the cluster category of an algebra of type $\mathbb{A}_3$ is a M\"obius strip. 
In fact the Auslander-Reiten quiver of the cluster category of any algebra of type $\mathbb{A}_n$ is a M\"obius strip for every $n \geq 2$.

One the one hand, they show that for any orientation of Dynkin quiver there is a one to one correspondence between cluster variables and indecomposable partial cluster-tilting objects; that there is a one to one correspondence between clusters and cluster-tilting objects; that the mutation is well-defined in all the indecomposable direct summands of any cluster-tilting object; and that the mutation of clusters and cluster-tilting objects are compatible. 
Note that these results were later generalised to the general case \cite{BCKMRT, CK}.

On the other, they showed that there is a natural inclusion of the module category of the path algebra into the cluster category such that every (partial) tilting module in the module category becomes a (partial) cluster-tilting object.
Moreover, they show that every possible mutation of tilting modules at the level of the module category becomes a mutation of cluster-tilting modules at the level of the cluster category. 

We said before that Happel and Ringel showed that much of the representation theory of tilted algebras can be described from the information we have about the representation theory of the hereditary algebras.
Now, the cluster categories associated to hereditary algebras have very nice properties, close to the properties of the hereditary algebras they come from. 
So Buan, Marsh and Reiten, emulating the construction of tilted algebras, introduced in \cite{BMR} the \textit{cluster-tilted algebras} as the endomorphism algebras of cluster-tilting objects in a cluster category. 
In this case, they showed that given a cluster-tilting object $T$ in $\mathcal{C}_A$, the functor $\Hom_{\C_A}(T, -)$ induces a equivalence of categories between $\mod\, (\End_{\C_A}(T))^{op}$ and the quotient of $\C_A$ by the ideal $\mathcal{I}(\tau T)$ of all the morphisms that factor through $\tau T$ the Auslander-Reiten translation of $T$ .

We have mentioned already that the module category of any hereditary algebra $A$ is naturally immersed in its cluster category $\C_A$. 
Moreover, if $T$ is a tilting object in $\mod\, A$, it turns out that $T$ becomes a cluster-tilting object in $\C_A$ when we apply the natural embedding. 
Then starting from $T$ we can construct a tilted algebra $\End_A(T)$ and a cluster-tilted algebra $\End_{\C_A}(T)$.
The relation between $\End_A(T)$ and $\End_{\C_A}(T)$ and their module categories was studied by Assem, Br\"ustle and Schiffler in a series of papers \cite{ABS1, ABS2, ABS3, ABS4}.
Firstly, they showed that one can recover $\End_{\C_A}(T)$ from $\End_A(T)$ via a process that they called \textit{relation extension} which bypasses the cluster category $\C_A$.
Moreover, they have shown that every cluster-tilted algebra is the relation extension of a tilted algebra. 
They also have characterised all the tilted algebras that have an isomorphic relation extension using particular structures that can be found in the Auslander-Reiten quivers of cluster-tilted algebras which are deeply related to the structures described by Liu \cite{LiuCrit} and Skowro\'nski \cite{Skocrit} for tilted algebras.

In order to start the construction of the cluster category, Buan, Marsh, Reineke, Reiten and Todorov assumed that the quiver in the initial seed of the algebra is acyclic. 
However, there is no reason why one should start with an acyclic quiver.
From a cluster perspective, any quiver is equally valid, so it was expected for a similar cluster category to exist regardless of the quiver we choose at the start. 
The first problem with the more general quivers arising in cluster theory is that they have cycles, so their path algebras are infinite dimensional. 
Then in order to use something close to tilting theory, we need to form the quotient of this path algebra by the correct ideal of relations. 
This problem was solved by Derksen, Weyman and Zelevinsky \cite{DWZ1, DWZ2} when they build ideals arising from certain of potentials associated to a quiver. 
They have shown that associated to each quiver there exists a special potential, that they called \textit{non-degenerate}, such that one can categorify the cluster algebra associated to the quiver using their \textit{decorated} representations. 
Moreover, they went further and showed that there exists a notion of mutation of non-degenerate potentials that is compatible with the cluster mutations of the quivers. 
We note that the notion of decorated representation gives rise to a rich theory that is closely related to that of $\tau$-tilting theory and it can be considered as a precursor of the latter.

Now that we have the correct algebras associated to the cyclic quivers in cluster theory we would like to have their corresponding cluster categories. 
To build these categories is not obvious. 
The main problem being that the construction of Buan, Marsh, Reineke, Reiten and Todorov uses heavily the structure of the derived category of the algebra and some key properties used in their construction fail when the algebra is not hereditary. 
This problem was overcome by Amiot in \cite{Amiot2009}, where she used the theory of Ginzburg dg-algebras developed by Keller and Yang in \cite{KY} to construct a cluster category which is compatible with the other notions of cluster categories existing to that moment. 

All the phenomena of cluster algebras and the close parallelism with tilting theory pointed to the existence of another extension of classical tilting theory where we would be always allowed to perform mutations, this time without extending the module category. 
For hereditary algebras, the construction of this theory was performed by Ingalls and Thomas in \cite{Ingalls2009}, where they introduced the so-called \textit{support tilting modules}. 
To explain the notion of support tilting module, let us come back to the limitations of classical tilting theory. 

As we did before consider $A$ to be the path algebra of the linearly oriented $\mathbb{A}_3$ quiver. 
Then $T = \rep{1\\2\\3}\oplus\rep{1}\oplus\rep{3}$ is a tilting module in $\mod\, A$. 

Ideally, given any choice $M$ of an indecomposable direct summand of $T$, we would like to construct a new tilting module whose indecomposable direct summands other than $M$ are the same as those of $T$.
In other words, we would like to replace each indecomposable direct summand of $T$ by another indecomposable in such a way that the resulting module is again tilting. 

The summand $\rep{1}$ is replaceable, since we can change it by $\rep{2\\3}$ to obtain $\rep{1\\2\\3}\oplus \rep{2\\3} \oplus \rep{3}$ which is tilting. 

We can also mutate at the summand $\rep{3}$, because it can be replaced by $\rep{1\\2}$ to obtain the tilting module $\rep{1\\2\\3}\oplus \rep{1\\2} \oplus \rep{1}$.

However, we cannot replace $\rep{1\\2\\3}$ by any other indecomposable module to obtain a new tilting module. 
This is a consequence of a classical result obtained independently by Assem \cite{AssemTorsTilt} and Smal\o~\cite{SmaloTorsionPairs},  which implies that every indecomposable projective-injective object in $\mod\, A$ is a direct summand of any tilting module in $\mod\, A$. 
In particular, $\rep{1\\2\\3}$ can not be replaced because it is a projective-injective module in $\mod\, A$.

The solution found by Ingalls and Thomas was to drop $\rep{1\\2\\3}$ from $T$ altogether to obtain $T'=\rep{1}\oplus\rep{3}$ which clearly is not tilting.
However, it is tilting on its \textit{support algebra}, which is constructed by taking a quotient of $A$ by the ideal generated by the idempotent included in the annihilator $\operatorname{ann} T'$ of $T'$.

More generally, they showed that for a hereditary algebra the mutation is always possible if we allow our tilting modules not to be supported over every vertex of the algebra. 

Now, for more general algebras this construction fails again. 
For instance, if we take the algebra $A$ to be the path algebra of the quiver 
$$
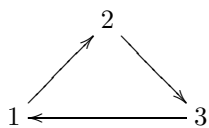
\xymatrix{
  & 2\ar[dr]& \\
  1\ar[ru] & & 3\ar[ll] }$$
modulo the ideal generated by all paths of length 2, we have that $A$ as a right module over itself is isomorphic to $\rep{1\\2}\oplus \rep{3\\1} \oplus \rep{2\\3}$.
Note that in this case every indecomposable projective is also injective. 
But at the same time we cannot drop any of the direct summands since the sum of the two remaining projective modules is supported on every vertex of the algebra. 

Something that we have not said before is that, by construction, in the cluster category we have that 
$$\Ext^1_{\C_A}(M,N) \cong \Hom_{\C_A}(N, \tau M).$$ 
This isomorphism can actually be translated to the module categories of non-hereditary cluster-tilted algebras. 
So, we can translate the cluster-tilting objects of the cluster category to the module category of a cluster-tilted algebra to get a series of modules which categorify perfectly the corresponding cluster algebra. 
However, these objects are not in general partial tilting objects because they might be of infinite projective dimension. 

Then, Adachi, Iyama and Reiten introduced $\tau$-tilting theory in \cite{AIR}, the object of study of these notes, by dropping the restriction on the projective dimension of the modules into consideration and replacing the classical rigidity with the notion of \textit{$\tau$-rigidity} that we will introduce in the next section. 
In doing so, as we will see in these notes, Adachi, Iyama and Reiten give a definition which can be easily checked in the module category of every finite dimensional algebra.
Moreover, as particular examples of this definition we can find the classical tilting modules and the modules over cluster-tilted algebras that we discussed in the previous paragraph. 

Before starting with the material of the lectures notes, we would like to point out that many results of $\tau$-tilting theory were developed independently by Derksen and Fei in \cite{DF2009}, where they studied \textit{general presentations} using methods of a more geometric nature.


\section{$\tau$-tilting theory: Basic definitions}\label{sec:definitions}

In this section we give the basic definitions of $\tau$-tilting theory. 
We also mention some of the basic relations between $\tau$-tilting theory and classical tilting theory. 
We start by giving the central definitions of this note: $\tau$-rigid and (support) $\tau$-tilting modules.

\begin{definition}\cite{Auslander1981, AuSmadd, AIR}\label{def:taurigid}
Let $A$ be an algebra and $M$ be an object in $\mod\, A$. 
We say that $M$ is $\tau$-rigid if $\Hom_A(M, \tau M) = 0$.
\end{definition}

\begin{definition}\label{def:supporttautilting}\cite{AIR}
Let $A$ be an algebra. 
A $\tau$-rigid $A$-module $M$ is $\tau$-tilting if $|M|=n$.
We say that a $\tau$-rigid $A$-module $M$ is support $\tau$-tilting if there exists an idempotent $e \in A$ such that $M$ is a $\tau$-tilting $A/AeA$-module, where $AeA$ is the two-sided ideal generated by $e$ in $A$.
\end{definition}

At first glance, $\tau$-tilting and tilting modules have little to do with each other.
If we compare Definition~\ref{def:tilting} with Definition~\ref{def:supporttautilting}, the only thing that a tilting and a $\tau$-tilting module have in common is that both are $A$-modules. 
However, there is a much deeper connection between the two concepts which follows from the so-called \textit{Auslander-Reiten formulas}.
Recall that \mbox{$D(-):= \Hom_k(-, k)$} denotes the classical duality functor, $\mathcal{I}_A(M, N)$ is the vector space of maps from $M$ to $N$ that factor through the injectives $A$-modules and $\mathcal{P}_A(M, N)$ is the vector space of maps from $M$ to $N$ that factor through the projectives $A$-modules.

\begin{theorem}\label{thm:AuslanderReitenformula}
Let $A$ be an algebra and let $M$ and $N$ be two $A$-modules.
Then there are funtorial isomorphisms
$$\Ext^1_A(M,N) \cong D\left(\frac{\Hom_A(N,\tau M)}{\mathcal{I}_A(N, \tau M)}\right) \cong D\left(\frac{\Hom_A(\tau^- N, M)}{\mathcal{P}_A(\tau^- N, M)}\right).$$
\end{theorem}

A module $M$ in $\mod\, A$ is said to be rigid if $\Ext^1_A(M,M)=0$. 
An immediate corollary of the Auslander-Reiten formulas is the following. 

\begin{corollary}
Let $M$ be an $A$-module. 
If $M$ is $\tau$-rigid then $M$ is rigid.
\end{corollary}

In fact this is the first of many other results relating tilting and $\tau$-tilting modules. 
In the following propositions we compile some properties relating the two notions.

\begin{proposition}\label{prop:properties1}
Let $A$ be an algebra and $T$ be a partial tilting module.
Then $T$ is $\tau$-rigid.
Moreover, if $T$ is tilting then $|T|=n$.
\end{proposition}

\begin{proposition}\label{prop:properties2}\cite{AsSS}
Let $M$ be a $\tau$-rigid module. Then the following hold.
\begin{enumerate}
	\item There are at most $n$ isomorphism classes of indecomposable direct summands of $M$. In short, $|M| \leq n$.
	\item If the annihilator $\operatorname{ann}(M)$ of $M$ is equal to the ideal $\{0\} \subset A$, then $M$ is a partial tilting module.
	\item If the projective dimension $pd\, M$ of $M$ is at most one, then $M$ is a partial tilting module. 
	\item If $|M|=n$ and $\operatorname{ann}(M) = \{0\}$, then $M$ is a tilting module.
\end{enumerate}
\end{proposition}

\begin{proposition}
Let $A$ be an algebra. 
Then an $A$-module $M$ is tilting if and only if $M$ is $\tau$-tilting and faithful.
\end{proposition}

\begin{proposition}
Let $A$ be an algebra. 
Then an $A$-module $M$ is tilting if and only if $M$ is $\tau$-tilting and $pd\, M \leq 1$.
\end{proposition}

The last result can be presented as evidence to the statement that $\tau$-tilting theory is a generalisation of tilting theory which is independent of the projective dimension of the objects.
Following this idea, in the last decade a series of works appeared generalising classical results in tilting theory to $\tau$-tilting theory.
Some of these results will be stated in Section~\ref{sec:tiltingtheorem}.

However, one needs to be careful when giving such statements, since some results on tilting theory do not hold in the context of $\tau$-tilting theory.
Also, before the definition of $\tau$-tilting theory there was at least one other generalisation of tilting theory to higher projective dimensions. 
We are referring to the generalised tilting modules introduced by Miyashita in \cite{Miyashita}.
They are defined as follows.

\begin{definition}\cite{Miyashita}\label{def:n-tilting}
Let $A$ be an algebra and $T$ be an $A$-module and $r$ be a positive integer. 
We say that $T$ is a $r$-tilting module if the following holds:
\begin{enumerate}
	\item[(i)] $\text{pd}_A T \leq r$, i.e. the projective dimension of $T$ is at most $r$.
	\item[(ii)] $\Ext^i_A(T,T)=0$ for all $1 \leq i \leq r$.
	\item[(iii)] There exists a short exact sequence of the form
	$$0 \to A \to T^{(1)} \to T^{(2)} \to \dots \to T^{(r)} \to 0$$
	where $T^{(i)}$ is a direct summand of a direct sum of copies of $T$ for all $1 \leq i \leq r$.
\end{enumerate}
\end{definition}

\medskip
We now give an example of all the support $\tau$-tilting modules in the module category of an algebra.

\begin{example}\label{ex:supandpair}
Let $A$ be the path algebra given by the quiver 
$$
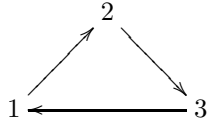
\xymatrix{
  & 2\ar[dr]& \\
  1\ar[ru] & & 3\ar[ll] }$$
modulo the second power of the ideal generated by all the arrows.
The Auslander-Reiten quiver of $A$ can be seen in Figure \ref{fig:Ar-quiverA3inclcon}.
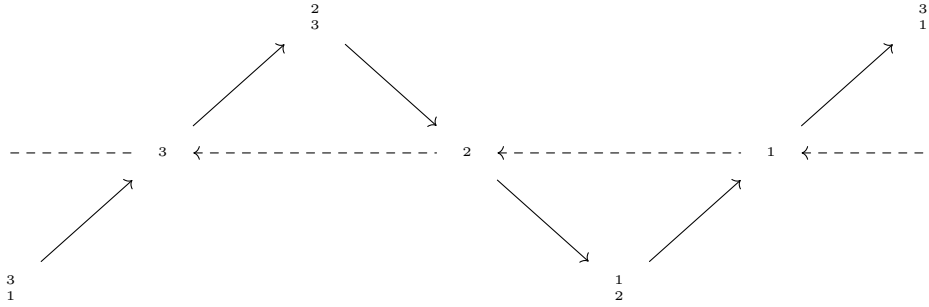
\begin{figure}[h]
    \centering
			\begin{tikzpicture}[line cap=round,line join=round ,x=2.0cm,y=1.8cm]
				\clip(-1.4,-0.3) rectangle (5.3,2.5);
					\draw [->] (-0.8,0.2) -- (-0.2,0.8);
					\draw [->] (3.2,0.2) -- (3.8,0.8);
					\draw [->] (0.2,1.2) -- (0.8,1.8);
					\draw [->] (4.2,1.2) -- (4.8,1.8);
					\draw [<-, dashed] (4.2,1.0) -- (5,1.0);
					\draw [dashed] (-1.,1.0) -- (-0.2,1.0);
					\draw [<-, dashed] (0.2,1.0) -- (1.8,1.0);
					\draw [<-, dashed] (2.2,1.0) -- (3.8,1.0);
					\draw [->] (2.2,0.8) -- (2.8,0.2);
					\draw [->] (1.2,1.8) -- (1.8,1.2);
				
				\begin{scriptsize}
					\draw[color=black] (-1,0) node {$\rep{3\\1}$};
					\draw[color=black] (5,2) node {$\rep{3\\1}$};
					\draw[color=black] (3,0) node {$\rep{1\\2}$};
					\draw[color=black] (4,1) node {$\rep{1}$};
					\draw[color=black] (0,1) node {$\rep{3}$};
					\draw[color=black] (2,1) node {$\rep{2}$};
					\draw[color=black] (1,2) node {$\rep{2\\3}$};
				\end{scriptsize}
			\end{tikzpicture}
\caption{The Auslander-Reiten quiver of $A$}
    \label{fig:Ar-quiverA3inclcon}

\end{figure}

Note that every module is represented by its Loewy series and both copies of $\rep{3\\1}$ should be identified, so the Auslander-Reiten quiver of $A$ has the shape of a Moëbius strip. 
In the first two columns of Table~\ref{table:supporttautilting} we give a complete list of the support $\tau$-tilting modules in $\mod\, A$ together with their associated idempotents.
\end{example}

Suppose that we are working on with the algebra of the previous example and we come across the module $M = \rep{1\\2} \oplus \rep{1}$.
After a quick calculation we can see that this module is not only $\tau$-rigid, but also support $\tau$-tilting with $e_3$ as its associated idempotent.
But at the same time, $M$ is a direct summand of the support $\tau$-tilting module $\rep{1\\2} \oplus \rep{1} \oplus \rep{3\\1}$. 
So, with this notation we cannot distinguish between the "complete" support $\tau$-tilting module $\rep{1\\2} \oplus \rep{1}$ and the "incomplete" $\tau$-rigid module $\rep{1\\2} \oplus \rep{1}$.

This can be solved using the notions of $\tau$-rigid and $\tau$-tilting pairs.
But before we give their definition, recall that given an idempotent $e \in A$ we have that the right ideal $eA$ is a projective module and that every projective arises this way. 

\begin{definition}
Let $A$ be an algebra, $M$ be an $A$-module and $P$ be a projective module. 
We say that the pair $(M,P)$ is $\tau$-rigid if $M$ is a $\tau$-rigid module and $\Hom_A(P,M)=0$.
A $\tau$-rigid pair is $\tau$-tilting if $|M|+|P|=n$.
\end{definition}

As you would expect, these two notations are equivalent. 
Indeed, given a support $\tau$-tiling module $M$ with associated idempotent $e$ we have that $(M, eA)$ is a $\tau$-tilting pair. 
Conversely, if $(M,P)$ is a $\tau$-tilting pair then we have that $P=eA$ for some idempotent $e\in A$. 
Then $M$ is a support $\tau$-tilting module with associated idempotent $e$.
The list of all $\tau$-tilting pairs of the algebra in Example~\ref{ex:supandpair} can be found in the third column of Table~\ref{table:supporttautilting}.

\begin{table}
\centering
\begin{tabular}{|c|c||c|}\hline
 support $\tau$-tilting module & idempotent & $\tau$-tilting pair \\\hline
 $\rep{1\\2} \oplus \rep{2\\3} \oplus \rep{3\\1}$ & $\varnothing$ &  $(\rep{1\\2} \oplus \rep{2\\3} \oplus \rep{3\\1}, 0)$ \\\hline
 $\rep{1\\2} \oplus \rep{2\\3} \oplus \rep{2}$ & $\varnothing$ & $(\rep{1\\2} \oplus \rep{2\\3} \oplus \rep{2}, 0)$\\\hline
 $\rep{1\\2} \oplus \rep{3\\1} \oplus \rep{1}$ & $\varnothing$ &  $(\rep{1\\2} \oplus \rep{3\\1} \oplus \rep{1}, 0)$  \\\hline
 $\rep{2\\3} \oplus \rep{3\\1} \oplus \rep{3}$ & $\varnothing$ & $(\rep{2\\3} \oplus \rep{3\\1} \oplus \rep{3}, 0)$ \\\hline
 $\rep{3\\1} \oplus \rep{3}$ & $e_2$ & $(\rep{3\\1} \oplus \rep{3}, \rep{2\\3})$ \\\hline
 $\rep{3\\1} \oplus \rep{1}$ & $e_2$ & $(\rep{3\\1} \oplus \rep{1}, \rep{2\\3})$ \\\hline
 $\rep{1\\2} \oplus \rep{1}$ & $e_3$ &  $(\rep{1\\2} \oplus \rep{1}, \rep{3\\1})$ \\\hline
 $\rep{1\\2} \oplus \rep{2}$ & $e_3$ &  $(\rep{2\\3} \oplus \rep{3}, \rep{1\\2})$ \\\hline
 $\rep{2\\3} \oplus \rep{3}$ & $e_1$ &  $(\rep{2\\3} \oplus \rep{3}, \rep{1\\2} )$ \\\hline
 $\rep{2\\3} \oplus \rep{2}$ & $e_1$ &  $(\rep{2\\3} \oplus \rep{2}, \rep{1\\2} )$ \\\hline
 $\rep{1}$ & $e_2 + e_3$ & $(\rep{1}, \rep{2\\3} \oplus \rep{3\\1})$\\\hline
 $\rep{2}$ & $e_1 + e_3$ & $(\rep{2}, \rep{1\\2} \oplus \rep{3\\1})$ \\\hline
 $\rep{3}$ & $e_1 + e_2$ & $(\rep{3}, \rep{1\\2} \oplus \rep{2\\3})$ \\\hline
 $\rep{0}$ & $e_1 + e_2 + e_3$ & $(\rep{0}, \rep{1\\2}\oplus \rep{2\\3} \oplus \rep{3\\1})$ \\\hline
\end{tabular}
\vspace{0.2cm}
\caption{Support $\tau$-tilting modules and $\tau$-tilting pairs in $\mod\, A$}
\label{table:supporttautilting}
\end{table}


\section{$\tau$-tilting pairs and torsion classes}\label{sec:torsionclasses}

In this section, after recalling the definition of torsion pairs and their basic properties, we will investigate the deep relation existing between $\tau$-tilting theory and torsion classes. 

\subsection{Torsion pairs and torsion classes}
The notion of torsion pairs, also known as torsion theories, started almost with the introduction of abelian categories as a generalisation of a well-known phenomenon in the category of finitely generated abelian groups, one of the most iconic examples of abelian categories. 

A classical classification result states that a finitely generated abelian group, up to isomorphism, has a unique torsion subgroup such that the resulting factor group is torsion-free. 
The extension of this fact to every abelian category was done by Dickson in \cite{Dickson1966} as follows.

\begin{definition}\label{def:torsionpair}\cite{Dickson1966}
Let $\mathcal{A}$ be an abelian category and let $(\T,\F)$ be a pair of subcategories of $\mathcal{A}$.
We say that $(\T,\F)$ is a torsion pair in $\A$ if the following holds.
\begin{enumerate}
	\item $\Hom_\A(X,Y)=0$ for all $X \in \T$ and $Y \in \F$.
	\item For all objects $M \in \A$ there exists, up to isomorphism, a unique short exact sequence 
	$$ 0 \to tM \to M \to fM \to 0$$
	such that $tM$ is an object of $\T$ and $fM$ is an object of $\F$.
\end{enumerate}
If $(\T, \F)$ is a torsion pair in $\A$ we say that $\T$ is a torsion class and $\F$ is a torsion-free class. 
Moreover, for each object $M$ of $\A$, we say that $$ 0 \to tM \to M \to fM \to 0$$ is the canonical short exact sequence of $M$ and that $tM$ is the torsion object of $M$ with respect to the torsion pair $(\T,\F)$.
\end{definition}

The previous definition is valid for an arbitrary abelian category. 
However, in these notes we are interested in a particular class of abelian categories, namely the categories of finitely generated modules over a finite dimensional algebra. 
These categories have many extra properties (for example they are length categories) that allow us to describe the torsion pairs more precisely.

\begin{proposition}
Let $A$ be an algebra. Then the following hold.
\begin{enumerate}
	\item A subcategory $\T$ of $\mod\, A$ is a torsion class if and only if $\T$ is closed under quotients and extensions. Moreover, in this case the torsion-free class associated to $\T$ is $$\F = \{Y \in \mod\, A : \Hom_A(X,Y)=0 \text{ for all $X \in \T$}\}.$$
	\item A subcategory $\F$ of $\mod\, A$ is a torsion-free class if and only if $\F$ is closed under subobjects and extensions. Moreover, in this case the torsion-free class associated to $\F$ is $$\T := \{X \in \mod\, A : \Hom_A(X,Y)=0 \text{ for all $Y \in \F$}\}.$$
\end{enumerate}
\end{proposition}

Suppose that $M$ is an $A$-module. 
Then we can ask the following: Is there a minimal torsion class in $\mod\, A$ containing $M$? 
The following result answers this question affirmatively. 

\begin{proposition}
Let $A$ be an algebra.
Then the intersection of arbitrarily many torsion classes is a torsion class. 
Likewise, the intersection of arbitrarily many torsion-free classes is a torsion-free class. 
\end{proposition}

Then, by the previous proposition, the minimal torsion class containing a given object $M$ of $\mod\, A$ is simply the intersection of all torsion classes containing $M$.  
Now, there is a more descriptive answer to this question, but to give that answer we need to introduce some notation. 

Let $\X$ be a subcategory of $\mod\, A$. 
The category $\Filt(\X)$ of objects filtered by $\X$ is defined as the category of all the objects $Y$ in $\mod\, A$ that admit a filtration 
$$0 = Y_0 \subset Y_1 \subset \dots \subset Y_{r-1} \subset Y_r = Y$$
such that the successive quotients $Y_{i}/ Y_{i-1}$ are objects in $\X$.
Note that $\Filt(\X)$ is the category of all the objects that can be constructed by making finitely many extensions by objects in $\X$. 
In other words, $\Filt(\X)$ is the extension closure of $\X$.

We define the category $\Fac \X$ as the category of objects $Y$ such that there exists an object $X$ in $\X$ and an epimorphism $p : X \to Y \to 0$.
Often in the notes, the category $\X$ we will be the additive category $add\, M$ additively generated by a module $M$.
In this case, by abuse of notation we will write $\Fac\, M$ instead of $\Fac (add\, M)$. 
Note that $\Fac\, M$ can be described as 
$$\Fac\, M = \{Y \in \mod\, A : \text{ there is an epimorphism }  p:M^r\to Y \to 0 \text{ for some $r \in \mathbb{N}$}\}.$$

Now we are able to give a better description of the minimal torsion class containing $M$.
As far as we know, the following result was a part of folklore but was first written down formally in \cite{DIJ}.

\begin{proposition}\label{prop:minimaltors}
Let $A$ be an algebra and $M$ be an $A$-module.
Then $\Filt (\Fac\, M)$ is the minimal torsion class containing $M$. 
\end{proposition}

\begin{remark}
Note that, in general, $\Fac(\Filt M)$ is \textbf{not} a torsion class since it might not be closed under extensions. 
\end{remark}

\subsection{Functorially finite torsion pairs and $\tau$-tilting theory}
From the previous subsection we have that to get the minimal torsion class containing $M$ one needs to first calculate $\Fac\, M$ and then take the extension closure of this category. 
However, sometimes $\Fac\, M$ is already closed under extensions, which makes the second step of the construction superfluous. 

The following theorem, originally proved by Auslander and Smal{\o} in \cite{AuSmadd}, is arguably the first result on $\tau$-tilting theory, even if this theory was formally introduced thirty years later. 

\begin{theorem}\label{thm:AuslanderSmalo1}
Let $A$ be an algebra and $M$ be an object in $\mod\, A$.
Then $\Fac\, M$ is a torsion class if and only if $M$ is $\tau$-rigid, that is $\Hom_A(M, \tau M)=0$.
Moreover, in this case 
$$M^\perp:= \{X \in \mod\, A : \Hom_A(M, X)=0 \}$$
is the torsion-free class such that $(\Fac\, M, M^\perp)$ is a torsion pair in $\mod\, A$.
\end{theorem}

As we just said, many years passed between the publication of this result and the start of $\tau$-tilting theory as an independent subject in representation theory. 
However, this was not the only result that worked with $\tau$-rigid objects. 
In fact, a well-established technique used in classical tilting theory to determine if an object is tilting was to show that the candidate $M$ was a $\tau$-rigid module such that $pd\, M \leq 1$ and $|M|=n$.
The interested reader is encouraged to surf the literature to look for such examples.

From a torsion theoretic point of view, the breakthrough made by Adachi, Iyama and Reiten in \cite{AIR} is that they showed that $\tau$-tilting pairs characterised a particular class of torsion classes, the \textit{functorially finite} torsion classes.

Let $\X$ be a subcategory of an abelian category $\A$ and suppose that $X$ is an object of $\X$ and $M$ is an arbitrary object of $\A$.
A morphism $f: X \to M$ in is called a \textit{right $\X$ -approximation} of M if any map $f':X' \to M$ with $X' \in \X$ factors through $f$. 
Dually, a morphism $g: M \to X$ is called a \textit{left $\X$-approximation} of M if any map $g':M \to X'$ with $X' \in \X$ factors through $g$.  
We say that $\X$ is \textit{contravariantly finite} (resp. \textit{covariantly finite}) if any object $M$ in $\A$ admits a right (resp. left) $\X$-approximation.
We say that $\X$ is functorially finite if it is both contravariantly finite and covariantly finite.

An important consequence of the uniqueness up to isomorphism of the canonical exact sequence of an object with respect to a torsion pair is the following.

\begin{proposition}
Let $(\T, \F)$ be a torsion pair in an abelian category $\A$ and let $M$ be an object of $\A$. If 
$$ 0 \to tM \to M \to fM \to 0$$
is the short exact sequence of $M$ with respect to $(\T, \F)$ then the canonical inclusion  $i : tM \to~M$ is a right $\T$-approximation and the canonical projection $p : M \to fM$ is a left $\F$-approximation.
In particular every torsion class in $\A$ is contravariantly finite and every torsion-free class $\F$ in $\A$ is covariantly finite.
\end{proposition}

Given a $\tau$-rigid module $M$, we know by Theorem~\ref{thm:AuslanderSmalo1} that $\Fac\, M$ is a torsion class. 
In fact, it turns out that $\Fac\, M$ is functorially finite, as shown by Auslander and Smal{\o} in \cite{AuSmadd}. 
Before we give the precise statement of the theorem, we will introduce some notation that will be useful in the rest of the paper.
Given a subcategory $\X$ of $\mod\, A$, we say that an object $M$ in $\X$ is \textit{$\Ext$-projective} if $\Ext^1_A(M,X)=0$ for every object $X$ in $\X$.
For every  functorially finite torsion class $\T$ of $\mod\, A$ we define $\P(\T)$ to be $\P(\T):=T_A^0 \oplus T_A^1$ where $f_A: A \to T_A^0$ is the minimal left $\T$-approximation of $A$ as an object of $\mod\, A$ and $T_A^1$ is the cokernel of $f_A$.

\begin{theorem}\cite{AuSmadd}
Let $\T$ be a functorially finite torsion class in $\mod\, A$. 
Then $\P(\T)$ is an $\Ext$-projective object in $\T$ such that $T$ is an object in $\add (\P(\T))$ for all $\Ext$-projective modules $T$ in $\T$.
Moreover $\T= \Fac\, \P(\T)$.
In particular $\P(\T)$ is a $\tau$-rigid $A$-module.
\end{theorem}

This last result implies that every functorially finite torsion class is generated by a $\tau$-rigid module. 
This defines a well defined map $\Phi: \tau\text{-rp-}A \to \text{ftors-}A$ from the set $\tau\text{-rp-}A$ of all $\tau$-rigid pairs to the set $\text{ftors-}A$ of functorially finite torsion classes.
The main contribution of \cite{AIR} to this problem is the proof of the fact that this map is a bijection if we consider restrict $\Phi$ to the set $\tau\text{-tp-}A \subset \tau\text{-rp-}A$ of $\tau$-tilting pairs.
The precise statement is the following.

\begin{theorem}\cite{AIR}\label{thm:bijectiontors}
Let $A$ be an algebra.
Then the map $\Phi : \tau\text{-tp-}A \to \text{ftors-}A$ defined by $$\Phi (M,P) = \Fac\, M$$ is a bijection.
Moreover, the inverse $\Phi^{-1} :  \text{ftors-}A \to \tau\text{-tp-}A$ is defined as
$$\Phi^{-1}(\T) = (\mathcal{P}(\T), {}^{\perp_P} \T)$$
where ${}^{\perp_p} \T$ is a basic additive generator of the category of projective modules $P$ such that $\Hom_A(P, T)=0$ for all $T \in \T$.
\end{theorem}


\section{A ($\tau$-)tilting theorem}\label{sec:tiltingtheorem}

We have mentioned in Section~\ref{sec:history} that the term \textit{tilting theory} was coined by Brenner and Butler, who showed in \cite{BB} what is now known as the tilting theorem. 
In this section we give a precise statement of the tilting theorem. 
Afterwards we show how this result can be generalised to $\tau$-tilting theory and mention some of its limits. 

\subsection{A tilting theorem}

We start by recalling the definition of a tilting module.

\begin{definition}\cite{BB}\label{thm:tilting}
Let $A$ be an algebra and $T$ be an $A$-module. 
We say that $T$ is a tilting module if the following holds:
\begin{enumerate}
	\item[(i)] $\text{pd}_A T \leq 1$, the projective dimension of $T$ is at most $1$.
	\item[(ii)] $T$ is rigid, that is $\Ext^1_A(T,T)=0$.
	\item[(iii)] There exists a short exact sequence of the form
	$$0 \to A \to T' \to T'' \to 0$$
	where $T', T''$ are direct summands of direct sums of $T$.
\end{enumerate}
\end{definition}

It follows from Proposition~\ref{prop:properties1} that every tilting module $T$ is $\tau$-tilting. 
Then we have that $T$ has a torsion pair associated to it, namely $(\Fac\, T, T^\perp)$.
Now, the torsion class $\Fac\, T$ can be characterised homologically as follows. 

\begin{proposition}
Let $T$ be a tilting module. 
Then $\Fac\, M = \left\{  X\in \mod\, A : \Ext_A^1(T,X)=0 \right\}$.
\end{proposition}

Now, for every object $M$ of $\mod\, A$ we have that $\End_A(M)$ is a finite dimensional algebra. 
In this case, $M$ has a natural structure of a left $\End_A(M)$-module structure. 
For the rest of the section we denote by $B:=\End_A(T)$ the endomorphism algebra of a tilting module $T$. 
The following proposition indicates the importance of tilting objects. 

\begin{proposition}\cite{BB}
Let $T$ be a tilting $A$-module. Then $T$ is tilting as a left $B$-module. 
Moreover, $T$ induces a torsion pair $(\X(T), \Y(T))$ in the category $\mod\, B$ of right $B$-modules where
$$\X(T):= \{X \in \mod\, B : X\otimes_{B} T =0\}$$
$$\Y(T):= \{Y \in \mod\, B : \operatorname{Tor}_1^B(Y,T) =0\}$$
\end{proposition}

We are now able to state the tilting theorem of Brenner and Butler. 

\begin{theorem}\cite{BB}
Let $A$ be an algebra, $T$ be a tilting $A$-module and $B=\End_A(T)$.
Then the following hold. 
\begin{enumerate}
	\item The algebra $\End_B^{op}(T)$ is isomorphic to $A$.
	\item The functor $\Hom_A(T, -): \Fac\, T \to \Y(T)$ is an equivalence of categories with quasi-inverse $-\otimes_B T : \Y(T) \to \Fac\, T$.
	\item The functor $\Ext^1_A(T, -):  T^\perp \to \X(T)$ is an equivalence of categories with quasi-inverse $\operatorname{Tor}_1^B(-, T) : \X(T) \to T^\perp$.
\end{enumerate}
\end{theorem}

\subsection{A $\tau$-tilting theorem}

We have said in various places that $\tau$-tilting theory can be seen as a generalisation of tilting theory. 
Hence, one would expect the existence of a generalisation of the tilting theorem to $\tau$-tilting theory. 
This was achieved in \cite{Ttaurodajas} building on the results of \cite{Jasso2015} that we discuss in the next section. 
Before stating the result, we need to recall some basic facts. 
The first key observation is that if $A$ and $C$ are algebras such that $C$ is a quotient of $A$, then $\mod\, C$ is a full subcategory of $\mod\, A$.
This implies immediately that $\End_A(M) \cong \End_C(M)$ for every $C$-module $M$.

In this section, we fix a $\tau$-tilting module $T$, we denote by $\operatorname{ann} T$ the annihilator of $T$ and by $C:=A / \operatorname{ann} T$ the quotient algebra of $A$ by $\operatorname{ann} T$.

\begin{proposition}
If $T$ is a $\tau$-tilting $A$-module, then $T$ is a tilting $C$-module. 
\end{proposition}
As a consequence of this proposition and Theorem~\ref{thm:tilting} we have the existence of a torsion pair $((\Fac\, T)_C, (T^\perp)_C)$ in $\mod\, C$ and a torsion pair $(\X(T), \Y(T))$ in $\mod\, B$ and equivalences of categories between them.
Now, we would like to compare the torsion pair $((\Fac\, T)_A, (T^\perp)_A)$  induced by $T$ in $\mod\, A$ with $(\X(T), \Y(T))$ in $\mod\, B$. 
The main ingredient to do that comes from the following two results.

\begin{proposition}\cite{AIR}
Let $T$ be a $\tau$-tilting module. 
Then $T^\perp = \Sub(\tau T)$.
\end{proposition}

\begin{proposition}
Let $C$ be a quotient algebra of $A$ and let $M$ be a $C$-module. 
Then the Auslander-Reiten translation $\tau_C M$ of $M$ in $\mod\, C$ is a submodule of the Auslander-Reiten translation $\tau_A M$ of $M$ in $\mod\, A$.
\end{proposition}
\noindent
Hence, if $T$ is a $\tau$-tilting $A$-module we have that the torsion pair $(\Fac\, T, T^\perp)$ in $\mod\, A$ coincides with the torsion pair $(\Fac\, T, T^\perp)$ in $\mod\, C$ if and only if $\tau_A M \cong \tau_C M$.
Then a $\tau$-tilting version of the tilting theorem of Brenner and Butler reads as follows. 

\begin{theorem}\cite{Ttaurodajas}
Let $A$ be an algebra, $T$ be a $\tau$-tilting $A$-module, $B=\End_A(T)$ and $C= A/\operatorname{ann} T$.
Then the following hold. 
\begin{enumerate}
	\item The algebra $\End_B^{op}(T)$ is isomorphic to $C$.
	\item The functors $\Hom_A(T, -): \Fac\, T \to \Y(T)$ is an equivalence of categories with quasi-inverse $-\otimes_B T : \Y(T) \to \Fac\, T$.
	\item The functors $\Ext^1_A(T, -):  T^\perp \to \X(T)$ is an equivalence of categories with quasi-inverse $\operatorname{Tor}_1^B(-, T) : \X(T) \to T^\perp$ if and only if $\tau_A T \cong \tau_C T$.
\end{enumerate}
\end{theorem}

\begin{remark}
We note that there examples of $\tau$-tilting pairs such that $\tau_A T \cong \tau_C T$, such as the so-called $\tau$-slices introduced in \cite{Ttaurodajas}.

We also note that there is another generalisation of the tilting theorem of Brenner and Butler for $2$-term silting objects proved by Buan and Zhou in \cite{BuanZhou2016}.
\end{remark}


\section{$\tau$-tilting reduction}

In the last section we explained the deep relationship between $\tau$-tilting pairs and torsion pairs. 
In this section we will consider the problem of completing a $\tau$-rigid pair, which we consider in two steps. 
In the first subsection we will show that there are two torsion pairs which are naturally associated to every $\tau$-rigid pair. 
In the second subsection we give a characterisation of all the completions of a $\tau$-rigid pair developed in \cite{Jasso2015}.
We finish this section by mentioning some bijections between the torsion classes of different categories.

\subsection{The Bongartz completion of a $\tau$-rigid pair}
The choice of taking $\tau$-rigid modules to develop $\tau$-tilting theory is arbitrary. 
In fact, using $\tau^{-1}$-rigid modules, that is, modules $N$ such that $\Hom_A(\tau^{-1}N, N)=0$, we can develop a completely dual $\tau^{-1}$-tilting theory. 
See \cite{AIR}.
In this case, dual of Theorem~\ref{thm:AuslanderSmalo1} reads as follows.

\begin{theorem}\cite{AuSmadd}
Let $A$ be an algebra and $N$ be an object in $\mod\, A$. 
Then the category 
$$\Sub N :=\{ Y \in \mod\, A : \text{ there is an monomorphism}  i: 0 \to Y  \to N^r \text{ for some $r \in \mathbb{N}$}\}$$
is a torsion-free class if and only in $N$ is $\tau^{-1}$-rigid.
Moreover, in this case 
$${}^\perp N:= \{X \in \mod\, A : \Hom_A(X,N)=0 \}$$
is the torsion class such that $({}^\perp N, \Sub N)$ is a torsion pair in $\mod\, A$.
\end{theorem}

Now, take a non-projective $\tau$-rigid module $M$. 
Then it is easy to see that $\tau M$ is $\tau^{-1}$-rigid. 
Indeed, 
$$\Hom_A(\tau^{-1} \tau M, \tau M)=\Hom_A(M, \tau M)=0.$$
Hence there are two torsion classes naturally associated to $M$, namely $\Fac\, M$ and ${}^\perp\tau M$. 
In the following theorem we give some results regarding the relation between these two torsion classes and $M$, all of which appeared already in \cite{AIR}.

\begin{theorem}\label{thm:torsmod}\cite{AIR}
Let $A$ be an algebra and $M$ be a $\tau$-rigid $A$-module.
Then the following holds.
\begin{enumerate}
	\item $\Fac\, M \subset {}^\perp \tau M $.
	\item The torsion classes $\Fac\, M$ and ${}^\perp \tau M$ coincide if and only if $M$ is $\tau$-tilting.
	\item Suppose that $\T$ is a functorially finite torsion class. Then $M$ is a direct summand of $\mathcal{P}(\T)$ if and only if $\Fac\, M \subset \T \subset {}^\perp \tau M$.
\end{enumerate}
\end{theorem}

From the previous theorem we have that ${}^\perp \tau M$ is the maximal torsion class having $M$ as an $\Ext$-projective, which makes the $\tau$-tilting module $\P({}^\perp \tau M)$ special enough to have a name.
We say that $\P({}^\perp \tau M)$ is the \textit{Bongartz completion} of $M$. 
This name was chosen because $\P({}^\perp \tau M)$ plays an analogous role in $\tau$-tilting theory to that of the Bongartz completion in the classical tilting theory.
To be more precise, Bongartz showed in \cite{BongartzComp} that if $T$ is a partial tilting module, then $\P({}^\perp \tau T)$ is a tilting module having $T$ as a direct summand. 
In other words, Bongartz showed that every partial tilting module can be completed to a tilting module.

If we use the language of $\tau$-rigid pairs instead of $\tau$-rigid modules we can be more precise in our statements. 
From now on, by abuse of notation, we say that a $\tau$-rigid pair $(M, P)$ is a direct summand of $(M',P')$ if $M$ is a direct summand of $M$ and $P$ is a direct summand of $P'$.

\begin{theorem}\cite{AIR}\label{thm:torspairs}
Let $A$ be an algebra and $(M,P)$ be a $\tau$-rigid pair in  $\mod\, A$.
Then the following hold.
\begin{enumerate}
	\item ${}^\perp \tau M \cap P^\perp$ is a torsion class and $\Fac\, M \subset {}^\perp \tau M  \cap P^\perp$.
	\item The torsion classes $M ^\perp$ and ${}^\perp \tau M  \cap P^\perp$ coincide if and only if $(M,P)$ is a $\tau$-tilting pair.
	\item Suppose that $\T$ is a functorially finite torsion class. Then $(M,P)$ is a direct summand of $\Phi^{-1}(\T)$ if and only if $\Fac\, M \subset \T \subset {}^\perp \tau M \cap P^\perp$.
\end{enumerate}
\end{theorem}

As for $\tau$-rigid modules, we say that $\Phi^{-1}({}^\perp \tau M \cap P^\perp)$ is the Bongartz completion of $(M,P)$.
But now we can also compute the $\tau$-tilting pair  $\Phi^{-1}(\Fac\, M)$ which is the $\tau$-tilting pair generating the smallest torsion class containing $M$.
In this case, we say that $\Phi^{-1}(\Fac\, M)$ is the Bongartz cocompletion of $(M,P)$.

\subsection{$\tau$-tilting reduction and torsion classes}

In this subsection we consider the problem of finding all $\tau$-tilting pairs having a given $\tau$-rigid pair $(M,P)$ as a direct summand. 
This problem was solved by Jasso in \cite{Jasso2015} using a procedure that is now known as \textit{$\tau$-tilting reduction}. 
Here we give a brief summary of that process. 

By Theorem~\ref{thm:torspairs} one knows that $(M,P)$ yields the torsion classes $\Fac\, M$ and $^\perp(\tau M) \cap P^\perp$. 
Moreover, Theorem~\ref{thm:torspairs} states the existence of a $\tau$-tilting pair of the form $(M\oplus M', P)$ such that $\Fac (M\oplus M')= {}^\perp(\tau M)\cap P^\perp$.

Now define $B_{(M,P)}= \End_A(M\oplus M')$ to be the endomorphism algebra of $M \oplus M'$.
In the algebra $B_{(M,P)}=\text{End}_A(M\oplus M')$, there is an idempotent element $e_{(M,P)}$ associated to the $B_{(M,P)}$-projective module $\Hom_{A}(M\oplus M', M)$. 
We define the algebra $\tilde{B}_{(M,P)}$ as the quotient of $B_{(M,P)}$ by the ideal generated by $e_{(M,P)}$, that is,
$$\tilde{B}_{(M,P)}:=B_{(M,P)}/B_{(M,P)} e_{(M,P)} B_{(M,P)}.$$

Now we are able to state one of the main results of \cite{Jasso2015}.

\begin{theorem}\cite{Jasso2015}\label{thm:tautiltred}
Let $(M,P)$ be a $\tau$-rigid pair in $\emph{mod} A$. 
Then the functor $$\emph{Hom}_A (M\oplus M', -): \emph{mod}\, A \to \emph{mod}\, B_{(M,P)}$$ induces an equivalence of categories 
$$F:M^\perp \cap {}^\perp\tau M \cap P^\perp \to \emph{mod}\, \tilde{B}_{(M,P)}$$ between the perpendicular category $M^\perp \cap {}^\perp\tau M \cap P^\perp$ of $(M,P)$ and the module category $\emph{mod}\, \tilde{B}_{(M,P)}$.
\end{theorem}

As a direct consequence of Theorem~\ref{thm:tautiltred} and Theorem~\ref{thm:torspairs} we obtain the following result. 

\begin{theorem}\label{thm:bijtorsJ}\cite{Jasso2015}
Let $(M,P)$ be a $\tau$-rigid pair in $\mod\, A$ and $\tilde{B}_{(M,P)}$ as above. 
Then the the functor 
$$\emph{Hom}_A (M\oplus M', -): \emph{mod}\, A \to \emph{mod}\, B_{(M,P)}$$ 
induces a bijection between the torsion classes $\T$ in $\mod\, A$ such that $\Fac\, M \subset \T \subset {}^\perp\tau M \cap P^\perp $ and the torsion classes in $\mod\, \tilde{B}_{(M,P)}$.

In particular the functor
$$\emph{Hom}_A (M\oplus M', -): \emph{mod}\, A \to \emph{mod}\, B_{(M,P)}$$ 
induces a bijection between the $\tau$-tilting pairs in $\mod\, A$ having $(M,P)$ as a direct summand and the $\tau$-tilting pairs in $\mod\, \tilde{B}_{(M,P)}$.
\end{theorem}

\begin{remark}
Note that Theorem~\ref{thm:bijtorsJ} does not give a specific number of completions of a given $\tau$-rigid pair $(M,P)$.
This is due to the fact that the number of $\tau$-tilting pairs in two algebras might differ hugely.
\end{remark}

\subsection{Bijections of torsion classes}
Note that the perpendicular category $M^\perp \cap {}^\perp\tau M \cap P^\perp$ defined by Jasso is a the intersection of the torsion class ${}^\perp\tau M \cap P^\perp$ with the torsion-free class $M^\perp$. 
Moreover, in this case $M^\perp\cap{}^\perp\tau M \cap P^\perp$ is what is called a \textit{wide subcategory} of $\mod\, A$. 
A subcategory $\X$ is called wide when it is closed under kernels, cokernels and extensions. 
In particular, this implies that $\X$ is an abelian category. 
Then Asai and Pfeiffer found in \cite{AsaiP2019} the following generalisation of Theorem~\ref{thm:bijtorsJ}.

\begin{theorem}\label{thm:bijtorsAP}\cite{AsaiP2019}
Let $(\T_1, \F_1)$ and $(\T_2, \F_2)$ be two torsion pairs in $\mod\, A$ such that $\T_1 \subset \T_2$. 
Suppose moreover that $\T_2 \cap \F_1$ is a wide subcategory of $\mod\, A$.
Then there is a bijection between the torsion classes $\T$ in $\mod\, A$ such that $\T_1 \subset \T \subset \T_2$ and the torsion classes in $\T_2 \cap \F_1$ given by map $\T \mapsto \T\cap\T_1$.
\end{theorem}
\smallskip

Note that, the intersection of a torsion class with a torsion-free class is not always a wide subcategory. 
However, it does always have some structure, namely that of a \textit{quasi-abelian} subcategory \cite{Tattar2020}.
The definition of quasi-abelian subcategories is a bit technical and it will be skipped. 

However, it is worth mentioning that Tattar showed in \cite{Tattar2020} that there is a well-defined notion of torsion classes in quasi-abelian subcategories. 
Moreover he showed that Theorem~\ref{thm:bijtorsJ} can be generalised to this setting as follows.

\begin{theorem}\label{thm:bijtorsT}\cite{Tattar2020}
Let $(\T_1, \F_1)$ and $(\T_2, \F_2)$ be two torsion pairs in $\mod\, A$ such that $\T_1 \subset \T_2$. 
Then there is a bijection between the torsion classes $\T$ in $\mod\, A$ such that $\T_1 \subset \T \subset \T_2$ and the torsion classes in $\T_2 \cap \F_1$ given by map $\T \mapsto \T\cap\T_1$.
\end{theorem}

\section{Torsion classes, wide subcategories and semibricks}\label{sec:bricks}

We have seen already a bijection between $\tau$-tilting pairs and functorially finite torsion classes. 
In this section we will see the relation of $\tau$-tilting theory with two other notions, namely wide subcategories and semibricks, which were described by Marks and Stovicek \cite{MarksStovicek} and Asai \cite{AsaiSemibricks}, respectively.

\subsection{Torsion classes and wide subcategories}

We start this subsection by recalling the definition of \textit{wide subcategories}\footnote{We give here the definition of wide subcategories as usually known in the context of representation theory of algebras. We warn the reader that in algebraic topology the same terminology has a different meaning.}. 

\begin{definition}
A subcategory $\W$ of $\mod\, A$ is said to be \textit{wide} if it is closed under extensions, kernels and cokernels. 
\end{definition}

The idea of this section is to describe maps between the set tors-$A$ of torsion classes in $\mod\, A$ and set wide-$A$ of wide subcategories in $\mod\, A$. 
Going from wide-$A$ to tors-$A$, the map that we take is quite natural. 
Indeed we define 
$$T(-) : \text{wide-}A \to \text{tors-}A$$
where $T(\W)$ is the minimal torsion class in $\mod\, A$ containing $\W$. 
We recall from Proposition~\ref{prop:minimaltors} that $T(\W)= \Filt(\Fac \W)$.
Actually this map has a really nice property. 
\begin{proposition}\cite{MarksStovicek}
The map $T(-) : \text{wide-}A \to \text{tors-}A$ is injective.
\end{proposition}

On the other direction we need to build a map $\alpha(-): \text{tors-}A \to \text{wide-}A$.
Given a torsion class $\T$ we define the subcategory $\alpha(\T)$ of $\mod\, A$ as follows. 
$$\alpha(\T):= \{ X\in \T : \ker f \in \T \text{ for all } f\in \Hom_A(Y,X) \text{ with $Y\in\T$} \}$$

\begin{proposition}\cite{MarksStovicek}
Let $\T$ be a torsion class in $\mod\, A$. 
Then $\alpha(\T)$ is a wide subcategory of $\mod\, A$.
\end{proposition}

Inside tors-$A$ there is the subset ftors-$A$ of functorially finite torsion classes, that we have already discussed before. 
Likewise, inside wide-$A$ there we have the subset of fwide-$A$ of wide subcategories $W$ of $\mod\, A$ such that which are functorially finite and such that $T(W)$ is a functorially finite torsion class.
The following theorem indicates what happens if we restrict the maps $\alpha(-)$  and $T(-)$ to these distinguished subsets. 

\begin{theorem}\cite{MarksStovicek}
The map $T(-) : \text{fwide-}A \to \text{ftors-}A$ is a bijection between the set \text{fwide-}A and $\text{ftors-}A$ with inverse $\alpha(-): \text{ftors-}A \to \text{fwide-}A$.
\end{theorem}

As an immediate consequence of the previous theorem and Theorem~\ref{thm:bijectiontors} we obtain the following corollary. 

\begin{corollary}
For every algebra $A$ there is a one-to-one correspondence between the set $\tau\text{-tp-}A$ of $\tau$-tilting pairs and the set $\text{fwide-}A$ of functorially finite wide subcategories of $\mod\, A$.
\end{corollary}

\subsection{Semibricks}

Let us recall the classical notion of bricks and the more novel notion of semibricks introduced by Asai in \cite{AsaiSemibricks}.

\begin{definition}
We say that an object $B$ in $\mod\, A$ is a \textit{brick} if its endomorphism algebra $\End_A(B)$ is a division ring. 
A set $\{B_1, \dots, B_t\}$ of bricks in $\mod\, A$ is said to be a semibrick if $\Hom_A(B_i, B_j)=0$ if $i$ is different from $j$.
\end{definition}

The classical example of a semibrick is the set $\{S(1), \dots, S(n)\}$ of non-isomorphic simple $A$ modules. 
Indeed, the classical Schur's lemma implies that $\{S(1), \dots, S(n)\}$ is a semibrick.
In fact, this is a particular case of a more general phenomenon. 
Recall that an object $M$ in a subcategory $\X$ of $\mod\, A$ is said to be \textit{relatively simple} if the only submodules of $M$ that are in $\X$ are $0$ and $M$.

Let $\W$ be a wide subcategory of $\mod\, A$ and let $B$ be a relative simple object in $\W$. 
Then it is easy to see that $S$ is necessarily a brick. 
Indeed, if $f\in \End_A(S)$ then im$f$ is a subobject of $B$ which is in $\W$ because $\W$ is wide. 
Hence $f$ is either the zero morphism or an isomorphism.
A similar argument shows that the set
$$\mathcal{S}(\W)=\{B: \text{ $B$ is relative simple in $\W$} \}.$$
is a semibrick.
On the other hand, one can see that given a semibrick $\cS$ the category $\Filt(\cS)$ is wide. 
In the following result $\text{sbrick-}A$ denotes the set of all semibricks in $\mod\, A$.

\begin{theorem}\cite{AsaiSemibricks}
The map $$\cS(-): \text{wide-}A \to \text{sbrick-}A$$ is a bijection with inverse $$\Filt(-): \text{sbrick-}A \to \text{wide-}A.$$
\end{theorem}

\noindent
Inside of $\text{sbrick-}A$ there are semibricks $\cS$ such that the minimal torsion class $T(\cS)$ containing $\cS$ is functorially finite.
We denote the set of all such semibricks by $\text{fsbrick-}A$. 
Combining the last theorem with the results of the previous subsection we obtain the following. 

\begin{theorem}\label{thm:somebijections}
There are bijections between the sets $\tau\text{-tp-}A$, $\text{ftors-}A$, $\text{fwide-}A$ and $\text{fsbrick-}A$.
\end{theorem}

There are several things that are worth mentioning here. 
Firstly, the bijections established in Theorem~\ref{thm:somebijections} are a very small subset of the existing bijections of interesting representation theoretic objects. 
A very nice paper that covered many of these bijections is the survey article \cite{BrustleYang} by Br\"ustle and Yang.
We choose to mention these here (and we will mention some other bijections later) because they were shown after the latest update of \cite{BrustleYang} and are not included there. 

We also want to emphasise that in the development of $\tau$-tilting theory a choice was made of working with the Auslander-Reiten translation $\tau$ and torsion classes in $\mod\, A$, instead of working with the inverse Auslander-Reiten translation $\tau^-$ and torsion free classes. 
However, the corresponding dual statements for these results hold.
In particular this adds many more bijections to this theory.

Finally, it is worth mentioning that Asai showed in \cite{AsaiSemibricks} an explicit bijection between $\text{fsbricks-}A$ and $\tau\text{-tp-}A$. 
This result was later recover in \cite{Treffinger2019} using the notions of $c$-vectors that will be discussed in a later section.

\section{The poset $\text{tors-}A$}\label{sec:poset}

The set $\text{tors-}A$ of all torsion classes in $\mod\, A$ has a natural poset structure, where the order is given by inclusion. 
In other words, given a pair of torsion classes $\T$ and $\T'$ we say that $\T \leq \T'$ if $\T \subset \T'$.
In this section we show some of the basic properties of this poset. 

\subsection{$\text{tors-}A$ is a complete lattice}
The aim of this subsection is to explain its title.
In order to do that we need to recall the definitions of a lattice and a complete lattice. 

\begin{definition}
A poset $\P$ is a \textit{lattice} if any two elements $x,y \in \P$ admit a greatest common lower bound $x\wedge y$, known as the \textit{meet} of $x$ and $y$, and a least common upper bound $x\vee y$, known as the \textit{join} of $x$ and $y$.
A lattice $\P$ is said to be \textit{complete} if every subset $S$ of $\P$ admits a lowest common greater bound $\bigvee_{x \in S} x$ (i.e. a \textit{a join}) and a greatest lower bound $\bigwedge_{x\in S} x$ (i.e. a \textit{meet}). 
\end{definition}

There has been a lot of work studying the lattice theoretic properties of the set $\text{tors-}A$ for an algebra $A$. 
One of the most important results in this direction was obtain by Demonet, Iyama, Reading, Reiten and Thomas in \cite{DIRRT}.

\begin{theorem}\cite{DIRRT}
The set $\text{tors-}A$ is a complete lattice for every algebra $A$. 
In this case, given a subset $S \subset \text{tors-}A$ the meet and the join of $S$ are defined as follows. 
$$\bigwedge_{\T\in S} \T := \bigcap_{\T \in S} \T$$
$$\bigvee_{\T \in S} \T:= T\left(\bigcup_{\T\in S} \T \right)$$
\end{theorem}

\begin{remark}
We note that the union of torsion classes is not always a torsion class since in general this union is not closed under extensions or quotients.
That is why we need to consider the minimal torsion class containing the union. 
\end{remark}

Now, suppose that $C$ is a quotient algebra of $A$. 
The following theorem tells us how $\text{tors-A}$ and $\text{tors-}C$ are related.

\begin{theorem}
If $A$ is an algebra and $C$ is a quotient of $A$, then there is an epimorphism of lattices $p: \text{tors-}A \to \text{tors-}C$. 
In other words, the lattice $\text{tors-}C$ is a quotient of the lattice $\text{tors-}A$.
\end{theorem}

One of the main ingredients of the proof of this theorem is the so-called \textit{brick labelling} that we will discuss in the next subsection.

\subsection{The Hasse quiver of $\text{tors-}A$}
We now shift the focus of our attention to the Hasse diagram of $\text{tors-}A$.
In order to do that, let us recall the notion of the Hasse quiver of a poset.
\begin{definition}
Given a poset $\P$, the Hasse quiver $\mathbf{H}(\P)$ of $\P$ is an oriented graph whose vertices correspond to the elements of $\P$ and there is an arrow $x \to y$ if $y \leq x$ such that $y\leq z \leq x$ implies that $x=z$ or $y=z$.
\end{definition}

In particular the arrows of the Hasse quiver $\mathbf{H}(\text{tors-}A)$ of $\text{tors-}A$ correspond to the \textit{maximally included} torsion classes. 
That is, for two torsion classes $\T$ and $\T'$  we say that $\T$ is maximally included in $\T'$ if $\T \subset \T'$ and $\T\subset \T'' \subset \T'$ implies that $\T''=\T$ or $\T''=\T'$.
It turns out that maximal inclusions of torsion classes have a very nice characterisation, as shown by Barnard, Carrol and Zhu in \cite{BCZ}. 

\begin{theorem}\cite{BCZ}\label{thm:minimalextendingmodule}
Let $\T$ and $\T'$ be two torsion classes such that $\T$ is maximally included in $\T'$. 
Then there exists a brick $B$ in $\mod\, A$ such that $\T' = \Filt( \T \cup \{B\})$.
In this case, we say that $B$ is the minimal extending module of the inclusion $\T \subset \T'$.
\end{theorem}

A direct consequence of this result is the following corollary, that it is usually known as the \textit{brick labelling} of $\mathbf{H}(\text{tors-}A)$.

\begin{corollary}
There is a well-defined labelling of the arrows of $\mathbf{H}(\text{tors-}A)$ by bricks in $\mod\, A$, where we label each arrow $\T' \to \T$ of $\mathbf{H}(\text{tors-}A)$ with the minimal extending module corresponding to the inclusion $\T \subset \T'$.
\end{corollary}

\begin{remark}
We note that brick labelling of $\mathbf{H}(\text{tors-}A)$ (or parts of it) can also be deduced as a consequence of several independent works that appeared simultaneously, namely \cite{AsaiSemibricks, BST2019, DIRRT, Treffinger2019}.
\end{remark}

We have that the set of functorially finite torsion classes $\text{ftors-}A$ is a subset of $\text{tors-}A$.
So the question is how their Hasse quivers $\mathbf{H}(\text{ftors-}A)$ and $\mathbf{H}(\text{tors-}A)$, respectively, are related. 
In order answer that question we need to know what happens when we have a maximal inclusion of torsion classes $\T \subset \T'$ such that either $\T$ or $\T'$ are functorially finite. 
A complete answer to this question is a consequence of the work of Demonet, Iyama and Jasso \cite{DIJ}.

\begin{theorem}\cite{DIJ}
Let $\T \subset \T'$ be a minimal inclusion of torsion classes in $\text{tors-}A$. 
Then $\T$ is functorially finite if and only if $\T'$ is functorially finite. 
\end{theorem}

A direct consequence of the previous theorem is the following.

\begin{corollary}
The Hasse quiver $\mathbf{H}(\text{ftors-}A)$ of $\text{ftors-}A$ and the Hasse quiver $\mathbf{H}(\text{tors-}A)$ of $\text{tors-}A$ are locally isomorphic.
\end{corollary}

\section{Mutations of $\tau$-tilting pairs and maximal green sequences}\label{sec:MGS}

As we said in in the introduction of these notes, $\tau$-tilting theory was conceived with the goal of completing the classical tilting theory with respect to mutation. 
In this section we start by discussing the notion of mutation from a torsion theoretic perspective. 
We finish it by speaking about the notion of maximal green sequences.
 
\subsection{Mutations of $\tau$-tilting pairs}
In order to speak about mutation we need to introduce a bit of notation. 

\begin{definition}
We say that a $\tau$-rigid pair $(M, P)$ is \textit{almost complete} if $|M|+|P|=n-1$. 
Given two $\tau$-tilting pairs $(M_1, P_1)$ and $(M_2, P_2)$, we say that $(M_1, P_1)$ is a mutation of $(M_2, P_2)$ if there is an almost complete $\tau$-rigid pair $(M,P)$ which is a direct summand of $(M_1,P_1)$ and $(M_2,P_2)$.
By abuse of notation we also say that $\Fac\, M_1$ is a mutation of $\Fac\, M_2$ if $(M_1, P_1)$ is a mutation of $(M_2, P_2)$.
\end{definition}

The main result about mutation of $\tau$-tilting pairs is the following theorem shown by Adachi, Iyama and Reiten in \cite{AIR}, that can be considered the most important result in $\tau$-tilting theory.

\begin{theorem}\cite{AIR}\label{thm:mutation}
Let $(M,P)$ be a $\tau$-tilting pair where $M= \bigoplus_{i=1}^k M_i$ and $P=\bigoplus_{j=k+1}^n P_j$ are basic modules. 

Then for every indecomposable direct summand $M_l$ of $M$ the almost complete $\tau$-rigid pair $(\bigoplus_{i\neq l} M_i, P)$ can be completed to a $\tau$-tilting pair different from $(M,P)$ in a unique way. 

Likewise, for every indecomposable direct summand $P_l$ of $P$ the almost complete $\tau$-rigid pair $(M, \bigoplus_{j\neq l} P_j)$ can be completed to a $\tau$-tilting pair different from $(M,P)$ in a unique way. 
\end{theorem}

\begin{remark}
We note that we can recover Theorem~\ref{thm:mutation} as a consequence of Theorem~\ref{thm:bijtorsJ}.
Indeed, one can verify that for all almost complete $\tau$-rigid pair $(M,P)$ the algebra $\tilde{B}_{(M,P)}$ is local, which implies that there is only one isomorphism class of simple modules in $\mod\, \tilde{B}_{(M,P)}$.
As a consequence of this, if $S$ is a simple module in $\mod\, \tilde{B}_{(M,P)}$ we have that $\Hom_{\tilde{B}_{(M,P)}}(X, S)=0$ implies that $X$ is isomorphic to zero.

Let $S$ be a simple module in $\mod\, \tilde{B}_{(M,P)}$ and let $\T$ be a torsion class in $\mod\, \tilde{B}_{(M,P)}$. 
Then we have two options: either $S\in \T$ or $S \not \in \T$.
If $S \in \T$, then $X \in \T$ for all $X \in \mod\, \tilde{B}_{(M,P)}$ since torsion classes are closed under extensions. 
Otherwise, we have that $\T = \{0\}$ by the argument above. 

This shows that every local algebra, no matter how complicated its representation theory, has exactly two torsion classes in its module category, which are the trivial torsion classes.
In particular, this implies that there are exactly two $\tau$-tilting pairs in $\mod\, \tilde{B}_{(M,P)}$ if $(M,P)$ is an almost complete $\tau$-rigid pair.
Hence Theorem~\ref{thm:mutation} follows from Theorem~\ref{thm:bijtorsJ}.
\end{remark}

Another important consequence of the previous argument is the following.

\begin{proposition}\cite{AIR}
Let $(M_1, P_1)$ and $(M_2, P_2)$ be $\tau$-tilting pairs that are mutations of each other and let $(M,P)$ be the almost complete $\tau$-rigid pair which is a direct summand of both.
Then either $\Fac\, M_1 \subset \Fac\, M_2$ or $\Fac\, M_2 \subset \Fac\, M_1$.
In particular, $\Fac\, M_1 \neq \Fac\, M_2$.
Moreover, if $\Fac\, M_1 \subset \Fac\, M_2$ then $\Fac\, M_1 = \Fac\, M$ and $\Fac\, M_2 = {}^\perp \tau M \cap P^\perp$.
\end{proposition}

\begin{remark}
From Theorem~\ref{thm:mutation} it follows that, given a $\tau$-tilting pair $(M,P)$ and a choice $(M'_1, P'_1)$ of an indecomposable direct summand of $(M,P)$, there is a unique $\tau$-tilting pair $(M_1,P_1)$ with the same indecomposable direct summands as $(M,P)$ except $(M’_1,P’_1)$. 
We say that $(M_1, P_1)$ is the mutation of $(M, P)$ at $(M'_1, P'_1)$.

Moreover, if $(M'_1, P'_1)$ and $(M'_2, P'_2)$ are distinct indecomposable direct summands of $(M,P)$, then the mutations $(M_1, P_1)$ and $(M_2, P_2)$ of $(M,P)$ at $(M'_1, P'_1)$ and $(M'_1, P'_1)$, respectively, are different.
In particular, this implies that for every $\tau$-tilting pair $(M,P)$ there are exactly $n$ $\tau$-tilting pairs which are a mutation of $(M,P)$ (see \cite{AIR}).
\end{remark}

In fact, there are explicit homological formulas to construct $(M_1, P_1)$ from $(M_2, P_2)$ and back but we will not explicit them here. 
The interested reader is encouraged to see \cite{AIR} for more details on the matter.

We now state a result by Demonet, Iyama and Jasso showing how the mutation of $\tau$-tilting pairs can be seen at the level of torsion theories. 

\begin{theorem}\cite{DIJ}\label{thm:chains}
Let $A$ be an algebra $(M,P)$ be a $\tau$-tilting pair and $\T$ be a torsion class in $\mod\, A$. Then the following hold.
\begin{enumerate}
	\item If $\T \subsetneq \Fac\, M$ then there exists a mutation $(M',P')$ such that $\T \subset \Fac\, M' \subsetneq \Fac\, M$.
	\item If $ \Fac\, M \subsetneq \T$ then there exists a mutation $(M'',P'')$ such that $\Fac\, M \subsetneq \Fac\, M'' \subset \T$.
\end{enumerate}
\end{theorem}

\smallskip
Using the notion of mutation of $\tau$-tilting pairs, one can construct a graph $\mathbf{\tau\text{-tp-}A}$ where the vertices are the $\tau$-tilting pairs in $\mod\, A$ and there is an edge between two $\tau$-tilting pairs if and only if they are mutation from each other. 
Using the results of this section and Theorem~\ref{thm:bijectiontors} we can prove the following result. 

\begin{proposition}
The graph $\mathbf{\tau\text{-tp-}A}$ is isomorphic to the undirected graph underlying $\mathbf{H}(\text{ftors-}A)$.
In particular, $\mathbf{H}(\text{ftors-}A)$ is an $n$-regular quiver.
\end{proposition}

\subsection{Maximal green sequences}
In the module category of any algebra $A$ there are always at least two torsion classes, sometimes called the trivial torsion classes, which are the whole of $\mod\, A$ and the torsion class $\{0\}$ containing only the objects that are isomorphic to the zero object.
Both of these torsion classes are functorially finite and it is not hard to see that $\Phi(A,0)= \mod\, A$ and $ \Phi(0,A)=\{0\}$.

Clearly $\{0\} \subsetneq \mod\, A$.
So we can apply Theorem~\ref{thm:chains}.1 and obtain a $\tau$-tilting pair $(M_1, P_1)$ which is a mutation of $(A, 0)$ such that $\{0\} \subset \Fac\, M_1 \subsetneq \mod\, A$.
If $\Fac\, M_1$ is not equal to $\{0\}$ we can repeat the process to obtain a mutation $(M_2, P_2)$ of $(M_1, P_1)$ such that $\{0\} \subset \Fac\, M_2 \subsetneq \Fac\, M_1 \subsetneq \mod\, A$.
We can repeat this process inductively to obtain a decreasing chain of torsion classes
$$\{0\} \subset \dots \subsetneq \Fac\, M_3 \subsetneq \Fac\, M_2 \subsetneq \Fac\, M_1 \subsetneq \mod\, A$$
which in general can continue forever. 
However, in some cases this process stops.
This leads to the following definition.

\begin{definition}
A \textit{maximal green sequence} of length $t$ in $\mod\, A$ is a finite set of $\tau$-tilting pairs $\{(M_i, P_i):  0 \leq i \leq t\}$ such that $(M_0, P_0)= (A,0)$, $(M_t, P_t)=(0,A)$ and $(M_i, P_i)$ is a mutation of $(M_{i-1}, P_{i-1})$ and $\Fac\, M_{i} \subset \Fac\, M_{i-1}$.
Equivalently, a maximal green sequence is a finite chain of torsion classes
$$\{0\} = \T_0  \subsetneq \T_1 \subsetneq \dots \subsetneq \T_{t-1} \subsetneq \T_t = \mod\, A$$
such that $\T_{i-1}$ is maximally included in $\T_i$ for all $1\leq i \leq t$.
\end{definition}

Maximal green sequences were originally introduced by Keller in \cite{Keller2011b} in the context of cluster algebras to give a combinatorial method to calculate certain geometric invariants known as Donaldson-Thomas invariants.
The interpretation of maximal green sequences in terms of chains of torsion classes was first used by Nagao in \cite{Nagao}. 
The definition given here can be considered as a generalisation to the setting of $\tau$-tilting theory of the original definition, since there are many examples of algebras which do not have a cluster counterpart, which first appeared in \cite{BST2019}.

We have preciously discussed in Section~\ref{sec:bricks} that given a maximal inclusion $\T \subset \T'$ of torsion classes $\T, \T'$ there is a brick $B$ in $\mod\, A$ such that $\T' = \Filt(\T \cup \{B\})$.
Then, applying this argument inductively we can see that if we have a maximal green sequence 
$$\{0\} = \T_0  \subsetneq \T_1 \subsetneq \dots \subsetneq \T_{t-1} \subsetneq \T_t = \mod\, A$$
we have a set of bricks $\mathcal{B}=\{B_1, B_2, \dots, B_t\}$ such that $\mod\, A = \Filt(\mathcal{B})$.
As it is, this is not a significant result since a small argument shows that all simple modules belong to $\mathcal{B}$ regardless of our starting maximal green sequence, and we know that simple modules filter every objet in $\mod\, A$.
However, one can show that the filtrations given in this way are unique in the following sense. 

\begin{theorem}
Let $\{0\} = \T_0  \subsetneq \T_1 \subsetneq \dots \subsetneq \T_{t-1} \subsetneq \T_t = \mod\, A$ be a maximal green sequence in $\mod\, A$ and let $\mathcal{B}=\{B_1, B_2, \dots, B_t\}$ be the set of bricks associated to it. 
Then for every non-zero object $M$ of $\mod\, A$ there is a filtration 
$$0 = M_0 \subset M_1 \subset \dots \subset M_{s-1} \subset M_s =M$$
such that $s \leq t$, $M_j / M_{j-1} \in \Filt(B_{j_i})$ and  $j_1 < j_2 < \dots < j_{s-1} < j_s$.
Moreover, this filtration is unique up-to-isomorphism.
\end{theorem}

\begin{remark}
The previous result is a consequence of two works that appeared independently, namely \cite{KellerSurvey} and \cite{T-HN-filt}.
In the appendix of \cite{KellerSurvey}, Demonet showed that the set $\mathcal{B}$ is what he calls a \textit{$I$-chain} and showed that every $I$-chain induces a unique filtration in every object of $\mod\, A$. 

A more general approach to this problem was given in \cite{T-HN-filt}, where it was shown that every chain of torsion classes induces a unique filtration for every object in $\mod\, A$. 
In that paper, this filtration was called the \textit{Harder-Narasimhan filtration} induced by the chain of torsion classes since it generalises the Harder-Narasimhan filtrations induced by stability conditions.
See also \cite{BST2022}.
\end{remark}

\begin{remark}
Note that the word \textit{green} in the name maximal green sequence does \textbf{not} make reference to any mathematician of name Green. 
Instead this word makes reference to the classical colouring in traffic lights.
The reason for this is that in cluster algebras there is no evident reason to say that a mutation is going forward or backwards. 
However, Keller \cite{Keller2011b} needed to impose such a direction to mutations in order to get the desired calculation. 
Then, he came up with a colouring of the vertices of the quiver associated to the cluster algebra in which a vertex is either green or red, which indicates if we are allowed to mutate at the given vertex or not, respectively. 
In this colouring, every vertex in the quiver of the initial seed is green and we are allowed to mutate at one green vertex at a time. 
The process finishes if after a finite number of mutations all the vertices in the quiver are red.
To learn more about this rich subject, see \cite{KellerSurvey}.
\end{remark}

\section{$\tau$-tilting finite algebras}\label{sec:tau-tiltingfinite}
To finish this part of the notes we speak about a new class of algebras that originated with the study of $\tau$-tilting theory, the so-called $\tau$-tilting finite algebras. 
They were introduced by Demonet, Iyama and Jasso as follows.

\begin{definition}\cite{DIJ}
An algebra $A$ is $\tau$-tilting finite if there are only finitely many $\tau$-tilting pairs in $\mod\, A$.
\end{definition}

Even though the class of $\tau$-tilting finite algebras has been recently introduced, they have received a lot of attention. 
In the following theorem we compile a series of characterisations of $\tau$-tilting finite algebras.

\begin{theorem}\label{thm:tautiltingfinite}
Let $A$ be an algebra. 
Then the following are equivalent.
\begin{enumerate}
	\item $A$ is $\tau$-tilting finite.
	\item There are finitely many indecomposable $\tau$-rigid objects in $\mod\, A$.
	\item \cite{DIJ} There are finitely many torsion classes in $\mod\, A$.
	\item \cite{DIRRT} There are finitely many bricks in $\mod\, A$.
	\item \cite{SchrollTreffinger} The lengths of bricks in $\mod\, A$ is bounded.
	\item \cite{AsaiSemibricks} There are finitely many semibricks in $\mod\, A$.
	\item \cite{MarksStovicek} There are finitely many wide subcategories in $\mod\, A$.
	\item \cite{Sentieri} Every brick in $\operatorname{Mod} A$ is a finitely presented $A$-module.
\end{enumerate}
\end{theorem}

\begin{remark}
Note that by Theorem~\ref{thm:tautiltingfinite}.3 we have that all torsion classes in the module category of a $\tau$-tilting finite algebra are functorially finite. 
\end{remark}
\smallskip

As we can see in Theorem~\ref{thm:tautiltingfinite}, $\tau$-tilting finite algebras have module categories that are somehow manageable from a torsion theoretic perspective, even if they are wild.  
As a consequence, there is an ongoing informal programme that aims to classify all the $\tau$-tilting finite algebras. 
This problem has been attacked by several people in different families of algebras. 
The following is a list of families of algebras where some progress to understanding on the problem has been made. 
Note that this list is not exhaustive nor efficient, since some families are included in others. 
\begin{itemize}
	\item \cite{Gabriel1972} Hereditary algebras.
	\item \cite{Adachi2016} Nakayama algebras.
	\item \cite{FZ2, BMR} Cluster-tilted algebras.
	\item \cite{Adachi2016a} Radical squared zero algebras.
	\item \cite{Zhang2017} Auslander algebras.
	\item \cite{Miz14} Preprojective algebras.
	\item \cite{PPP2018, Plamondon2019} Gentle algebras.
	\item \cite{AAC18} Brauer graph algebras.
	\item \cite{STV21} Special biserial algebras.
	\item \cite{August2020} Contraction algebras.
	\item \cite{Mousavand2019a} Non-distributive algebras.
\end{itemize}

\section{Brauer-Thrall conjectures and $\tau$-tilting theory}

The systematic study of the representation theory of finite dimensional algebras, one can argue, started around the 1940s. 
At that time it was already known that every $A$-module can be written as a direct summand of indecomposable $A$-modules in essentially one way. 
As a consequence, people started to classify the algebras in two types, \textit{representation finite} and \textit{representation infinite}, depending on their number of isomorphism classes of indecomposable modules.
Much of the motivation in the early days of representation theory of finite dimensional algebras stems from this quest of determining algebras of finite representation type with an important role being played by the first and second Brauer-Thrall Conjectures, which were proved subsequently by Ro\u{\i}ter \cite{Roiter1968} and Auslander \cite{Auslander1974} for the first and by Nazarova and Ro\u{\i}ter \cite{Nazarova1971} and Bautista \cite{Bautista1985} for second.
For an historical survey on Brauer-Thrall conjectures and its influence in representation theory, please see \cite{Gustafson1982}.
Their statement is the following.

\begin{conjecture}[First Brauer-Thrall Conjecture]\label{conj:BT1}\cite{Auslander1974, Roiter1968}
Let $A$ be an algebra. 
Then $A$ is of finite representation type if and only if there is a positive integer $d$ such that $\dim_\k (M) \leq d$ for every indecomposable $A$-module $M$.
\end{conjecture}

\begin{conjecture}[Second Brauer-Thrall Conjecture]\label{conj:BT2}\cite{Bautista1985, Nazarova1971}
If $A$ is an algebra of infinite representation type, then there is an infinite family of positive integers $\{d_i : i \in \mathbb{N} \}$ such that for every $d_i$ there is an infinite family of indecomposable $A$-modules $\{M^{d_i}_j\}$ where $\dim_\k (M_j^{d_i})=d_i$ for all $j$.
\end{conjecture}

Then using the characterisation of $\tau$-tilting finite algebras given in Theorem~\ref{thm:tautiltingfinite}.4 one can state a $\tau$-tilting analogue of the first Brauer-Thrall Conjecture by restricting the universe of modules to consider from indecomposable to bricks.
The statement of the conjecture, which was proven for every finite dimensional algebra over a field in \cite{SchrollTreffinger}, is the following.

\begin{theorem}[First $\tau$-Brauer-Thrall Conjecture]\label{conj:tBT1}\cite{SchrollTreffinger}
Let $A$ be an algebra. 
Then $A$ is $\tau$-tilting finite if and only if there is a positive integer $d$ such that $\dim_\k (M) \leq d$ for every brick $M$ in $\mod\, A$.
\end{theorem}
The statement of this conjecture appeared independently in \cite{STV21} and \cite{Mousavand2019}.
We note that recently a new proof of the validity of this conjecture was given in \cite{MP2021}.

\smallskip
One can see immediately that the $\tau$-tilting version of the first Brauer-Thrall conjecture is a direct translation of the original. 
However, one can not do the same in the case of the second Brauer-Thrall conjecture, since there are $\tau$-tilting infinite algebras having finitely many infinite families of bricks. 
The $\tau$-tilting version of the second Brauer-Thrall conjecture is the following.

\begin{conjecture}[Second $\tau$-Brauer-Thrall Conjecture]\label{conj:tBT2}
If $A$ is a $\tau$-tilting finite algebra, then there is a positive integer $d$ and an infinite family of bricks $\{M^{d}_j\}$ in $\mod\, A$ such that \mbox{$\dim_\k (M_j^{d})=d$} for all $j$.
\end{conjecture}

This conjecture is still open in general. 
However, it has been verified for gentle algebras \cite{Plamondon2019}, special biserial algebras in \cite{STV21, Mousavand2019} and for distributive algebras in \cite{Mousavand2019a}.

\begin{remark}
We note that if an algebra satisfies the second $\tau$-Brauer-Thrall conjecture then it also satisfies the second Brauer-Thrall conjecture by a result of Smal\o~\cite{Smalo1980}.
\end{remark}

\begin{remark}
We also note that there are many examples of $\tau$-tilting finite algebras that are of infinite representation type, such as preprojective algebras \cite{Miz14} or contraction algebras \cite{August2020}, to name just a few.
\end{remark}


\section{Dimension vectors and $g$-vectors}\label{sec:vectors}

In the introduction of these notes we said that many developments that occurred in representation theory in the twenty-first century, including $\tau$-tilting theory, were aiming to \textit{categorify} cluster algebras to some extent. 

In loose terms, the term \textit{categorification} refers to the process of explaining some combinatorial phenomena by showing the existence of some underlying categorical phenomena. 
For instance, the bijection between the indecomposable $\tau$-rigid modules in the module category of a hereditary algebra of Dynkin type and the non-initial variables in the cluster algebra of the corresponding Dynkin type is a categorification of cluster variables. 

But as there is a process of categorification, there is also a process of \textit{decategorification}, a process where you start with a category and you find some combinatorial or numerical data that reflects the phenomena occurring at the categorical level. 

From now on, we shift our focus to a different decategorification of $\tau$-tilting theory of an algebra using integer vectors.

\subsection{The Grothendieck group of an algebra}

The most classical decategorification using integer vectors of the representation theory associated to an algebra is the \textit{Grothendieck group} of a category. 
We start by recalling the definition in the case of arbitrary abelian categories.

\begin{definition}\label{def:Grothendieckgroup}
Let $\A$ be an abelian category. 
The Grothendieck group $K_0(\A)$ of $\A$ is the quotient of the free abelian group generated by the isomorphism classes $[M]$ of all objects $M \in \A$ modulo the ideal generated by the short exact sequences as follows. 
$$K_0(\A) = \frac{\langle[M] : M \in \A\rangle}{\langle[M]-[L]-[N]: 0 \to L \to M \to N \to 0 \text{ is a short exact sequence in $\A$}\rangle}$$
\end{definition}

In these notes we are interested only in the module categories $\mod\, A$ of finite dimensional algebras $A$ over an algebraically closed field. 
By abuse of notation, the Grothendieck group $K_0(\mod\, A)$ of $\mod\, A$ will be denoted by $K_0(A)$ and we will refer to it as the Grothendieck group of $A$.
As an immediate consequence of the Jordan-H\"older theorem for module categories we have the following result. 

\begin{theorem}
Let $A$ be an algebra. 
Then $K_0(A)$ is isomorphic to $\mathbb{Z}^n$, where $n$ is the number of isomorphism classes of simple modules in $\mod\, A$.
\end{theorem}

From now on, we fix a complete set $\{[S(1)], \dots, [S(n)]\}$ of isomorphism classes of simple $A$-modules. 
Clearly, $\{[S(1)], \dots, [S(n)]\}$ forms a basis of $K_0(A)$.
However, this is not the only basis of $K_0(A)$.
In these notes, when we speak about the Grothendieck group of $A$ we always assume that the basis chosen to represent our vectors is the basis given by the simple modules with a fixed order. 

\begin{theorem}
Let $A$ be an algebra and $K_0(A)$ be its Grothendieck group having as canonical basis the set $\{[S(1)], \dots, [S(n)]\}$ of isomorphism classes of simple $A$-modules. 
Then for every object $M \in \mod\, A$ we have that 
$$[M] = (\dim_\k(\Hom_A(P(1), M)), \dots, \dim_\k(\Hom_A(P(n), M)))$$
$$[M] = (\dim_\k(\Hom_A(M, I(1)), \dots, \dim_\k(\Hom_A(M, I(n)))))$$
where $P(i)$ and $I(i)$ are the projective cover and the injective envelope of the simple $S(i)$, respectively, for all $1 \leq i \leq n$.
\end{theorem}

The previous result justifies that the element of the Grothendieck group $[M]$ associated to $M$ is often called the \textit{dimension vector} of $M$, terminology that we adopt in these notes as well.

Sometimes in the literature one finds the notation \underline{dim}$M$ for the dimension vector, reserving $[M]$ for the abstract class of $M$ in the Grothendieck group with no prefered basis of $K_0(A)$.

\begin{remark}
In the previous result we are actually using the hypothesis that $A$ is an algebra over an algebraically closed field. 
Otherwise, the result is not true in general. 
We warn the reader that this remark is also valid for several other results in this section. 
\end{remark}

\subsection{g-vectors}

Another set of integer vectors that can be associated to the category of finitely presented $A$-modules is the set of $g$-vectors.

Although the idea of $g$-vectors has been around for several decades, their systematic study is rather recent since the main motivation behind their study, as the reader might be guessing already, lies in the categorification of cluster algebras. 

In fact, the name $g$-vector itself comes from cluster theory. 
The $g$-vectors were introduced by Fomin and Zelevinsky in \cite{FZ4}, where they conjectured that cluster variables could be parametrised using $g$-vectors. 
Later on, it was shown that $g$-vectors encoded the projective presentation of $\tau$-rigid $A$-modules. 
Their definition is the following. 

\begin{definition}
Let $M$ be an $A$-module. 
Choose the minimal projective presentation $$P_1\longrightarrow P_0\longrightarrow M\longrightarrow 0$$  of $M$, where $P_0=\bigoplus\limits_{i=1}^n P(i)^{a_i}$ and $P_1=\bigoplus\limits_{i=1}^n P(i)^{b_i}$. 
Then the $g$-vector of $M$ is defined as
$$g^M=(a_1-b_1, a_2-b_2,\dots, a_n-b_n).$$
The $g$-vector of a $\tau$-rigid pair $(M,P)$ is defined as $ g^{M}-g^{P}$.
\end{definition}

\begin{remark}\label{rmk:g-vect}
Recently, Nakaoka and Palu introduced in \cite{NakaokaPalu} introduced the notion of \textit{extriangulated categories}.
These categories are a generalisation of abelian where the notion of \textit{conflations} replace and generalise the notion of short exact sequences. 
So, one can define the Grothendieck group of any extriangulated category as we did in Definition~\ref{def:Grothendieckgroup}, replacing the short exact sequences by the conflations in this category. 
In the following section we will see that $g$-vectors can be thought as the elements of the Grothendieck group of an extriangulated category denoted by $K^{[-1,0]}(\proj\, A)$.
\end{remark}


\section{$\tau$-tilting theory and $2$-term silting complexes}\label{sec:silting}

Associated to any finite dimensional algebra $A$ there is a triangulated category known as the homotopy category of bounded complexes of finitely generated projective $A$-modules usually denoted by $K^{b}(\proj\, A)$. 
The theory of homotopy categories is very rich, but out of the scope of these notes. 
The interested reader is encouraged to see \cite{Miller2020} for a detailed account on that matter. 

In this section we first introduce an abstract category associated to every algebra $A$ that we call $K^{[-1,0]}(\proj\, A)$ because it can be identified with a full subcategory of $K^{b}(\proj\, A)$.
Later on, we explain the relation between  $K^{[-1,0]}(\proj\, A)$ and $\tau$-tilting theory.

\subsection{The category $K^{[-1,0]}(\proj\, A)$.}

As one usually does when defining a category, let us start by defining its objects.
In this case, an object $\PP$ of $K^{[-1,0]}(\proj\, A)$ is a complex of the form $\PP:= P_{-1} \xrightarrow{f} P_0$ where $P_0$ and $P_{-1}$ are two projective $A$-modules and $f \in \Hom_A(P_{-1}, P_0)$.
Now, given two objects $\PP:= P_{-1} \xrightarrow{f} P_0$ and $\QQ:= Q_{-1} \xrightarrow{f'} Q_0$ in $K^{[-1,0]}(\proj\, A)$ we define $g \in \Hom_{K^{[-1,0]}(\proj\, A)}(\PP, \QQ)$ as a pair $g:=(g_{-1}, g_0)$ such that $g_i \in \Hom_A(P_i, Q_i)$ making the following diagram commutative
$$\xymatrix{ P_{-1} \ar[r]^f\ar[d]_{g_{-1}} & P_0\ar[d]^{g_0}\ar@{.>}[dl]_{h} \\
		Q_{-1} \ar[r]^{f'} & Q_0}$$
modulo the equivalence relation given by if $g, g' \in \Hom_{K^{[-1,0]}(\proj\, A)}(\PP, \QQ)$, we impose that $g$ is equal to $g'$ if there exists a map $h: P_0 \to Q_{-1}$ such that $f'h=g_0$ and $hf=g_{-1}$.

We note that in this category one can associate to every pair of objects $\PP, \QQ \in K^{[-1,0]}(\proj\, A)$ an abelian group $\mathcal{E}_{K^{[-1,0]}(\proj\, A)}(\PP, \QQ)$, known as the group of \textit{conflations} of $\PP$ by $\QQ$.
This group is defined as the group of maps $t\in \Hom_A(P_{-1}, Q_0)$ \textit{up to homotopy}, meaning that $t, t'\in \Hom_A(P_{-1}, Q_0)$ coincide in this category if there are maps $t_i\in \Hom_A(P_i, Q_i)$ such that $t-t'=f'h_1 - h_0f$.
$$\xymatrix{ &P_{-1} \ar[r]^f \ar[d]^t\ar@{.>}[dl]_{h_1} & P_0\ar@{.>}[dl]^{h_0} \\
		Q_{-1} \ar[r]^{f'}& Q_0 &}$$

The notion of a conflation in an extriangulated category is a generalisation of that of short exact sequences in abelian categories and distinguished triangles in triangulated categories. 
As such, a conflation $t \in \mathcal{E}_{K^{[-1,0]}(\proj\, A)}(\PP, \QQ)$  can be realised as 
$$\mathbf{Q}\overset{a}\rightarrowtail \mathbf{E_t} \overset{b}\twoheadrightarrow \mathbf{P}$$
where $a \in \Hom_{K^{[-1,0]}(\proj\, A)}(\QQ, \mathbf{E_t})$, $b\in \Hom_{K^{[-1,0]}(\proj\, A)}(\mathbf{E_t}, \PP)$ and $\mathbf{E_t}$ is in $K^{[-1,0]}(\proj\, A)$.
The existence of conflations and their realisation endow the category $K^{[-1,0]}(\proj\, A)$ with the structure of a \textit{extriangulated category}, a notion very recently introduced by Nakaoka and Palu in \cite{NakaokaPalu}. 
This follows from the fact that $K^{[-1,0]}(\proj A)$ is equivalent to an extension-closed subcategory of the triangulated category $K^{b}(\proj A)$. 
See \cite{PPPP} for more details.

One can then define the notion of a projective in $K^{[-1,0]}(\proj\, A)$ as an object $\PP\in K^{[-1,0]}(\proj\, A)$ such that $\mathcal{E}_{K^{[-1,0]}(\proj\, A)}(\PP, \QQ)=0$ for every object $\QQ \in K^{[-1,0]}(\proj\, A)$.
Likewise, an injective object in $K^{[-1,0]}(\proj\, A)$ is an object $\mathbf{I}\in K^{[-1,0]}(\proj\, A)$ such that $\mathcal{E}_{K^{[-1,0]}(\proj\, A)}(\QQ, \mathbf{I})=0$ for every object $\QQ \in K^{[-1,0]}(\proj\, A)$.
As it turns out, the projective and injective objects of $K^{[-1,0]}(\proj\, A)$ can be described explicitly as follows.
$$\proj_{K^{[-1,0]}(\proj\, A)} = \left\{0 \xrightarrow{0} P : \text{ $P \in \proj\, A$} \right\} $$
$$\inj_{K^{[-1,0]}(\proj\, A)} = \left\{P \xrightarrow{0} 0  : \text{ $P \in \proj\, A$} \right\} $$

As we said before, we want to relate the category $K^{[-1,0]}(\proj\, A)$ with the $\tau$-tilting theory of $A$. 
To do that, we first need to know how one can go from $\mod\, A$ to $K^{[-1,0]}(\proj\, A)$ and back.

In fact, one can see $\mod\, A$ inside $K^{[-1,0]}(\proj\, A)$.
This comes from the fact that, by definition, $\mod\, A$ is the category of finitely presented $A$-modules. 
This means that every $A$-module $M$ admits a minimal projective presentation
$P_{-1} \xrightarrow{f} P_0 \xrightarrow{} M \xrightarrow{} 0$, which induces a map $\PP(-): \mod\, A \to K^{[-1,0]}(\proj\, A)$ defined  as $\PP(M)= P_{-1} \xrightarrow{f} P_0$ on the objects.
We note that if $A$ is a hereditary algebra, then $\PP(-)$ is actually a functor. 

On the other direction, there is a natural functor from $H_0(-): K^{[-1,0]}(\proj\, A) \to \mod\, A$, usually known as the \textit{0-th homology}, which in this case consists of $H_0(\PP)= \coker f$, where $\PP = P_{-1} \xrightarrow{f} P_0$.
In fact, one can show that the kernel of the functor $H_0$ corresponds to the ideal $\left<A \xrightarrow{} 0\right>$ generated by the injective objects of $ K^{[-1,0]}(\proj\, A)$. 

$$\mod\, A \cong  \frac{K^{[-1,0]}(\proj\, A)}{\left<A \xrightarrow{} 0 \right>}$$

It was noted by Gorsky, Nakaoka and Palu \cite{GorskyNakaokaPalu} that for every finite dimensional algebra $A$, the category $K^{[-1,0]}(\proj\, A)$ is hereditary. 
In other words, they have shown that for every object $\mathbf{M} \in K^{[-1,0]}(\proj\, A)$ there is a conflation of the form 
$$\PP_{1} \rightarrowtail \PP_0 \twoheadrightarrow \mathbf{M}$$
where $\PP_0$ and $\PP_1$ are projectives in $K^{[-1,0]}(\proj\, A)$.
This fact is highly surprising, since the global dimension of the module category of an algebra is arbitrary. 
Moreover, this and the results that we state in the following subsection give a partial explanation of the good behaviour of $\tau$-tilting theory for every algebra. 

It follows from the results in this section that the Grothendieck group $K_0(K^{[-1,0]}(\proj\, A))$ of $K^{[-1,0]}(\proj\, A)$ is isomorphic to $\mathbb{Z}^n$, where the set $[\PP(i)]=[0 \to P(i)]$ form a natural basis this Grothendieck group. 
In particular we have that the $g$-vector $g^M$ of any $M \in \mod\, A$ is the class $[\PP(M)]$ of $\PP(M)\in K^{[-1,0]}(\proj\, A)$ written in the basis $\{[\PP(1)], \dots, [\PP(n)]\}$.

\subsection{$\tau$-tilting theory in $K^{[-1,0]}(\proj\, A)$}

As we pointed out before, the objects of $\mod\, A$ can be seen as objects of $K^{[-1,0]}(\proj\, A)$.
In this subsection we see how the $\tau$-tilting theory of an algebra $A$ can be studied $K^{[-1,0]}(\proj\, A)$. 
For this we start with the following definition.  

\begin{definition}
Let $A$ be an algebra and $\PP$ be an object in $K^{[-1,0]}(\proj\, A)$.
We say that $\PP$ is \textit{presilting} if $\mathcal{E}_{K^{[-1,0]}(\proj\, A)}(\PP, \PP)=0$.
Moreover, we say that $\PP$ is \textit{silting} if the number of isomorphism classes of indecomposable direct summands of $\PP$ is equal to $|A|$.
\end{definition}

\begin{remark}
We said before that $K^{[-1,0]}(\proj\, A)$ can be identified with a full subcategory of $K^{[b}(\proj\, A)$. 
Under this identification, the previous definition coincides with a characterisation of $2$-term silting objects in $K^{b}(\proj\, A)$ given by Adachi, Iyama and Reiten in \cite{AIR}.
\end{remark}

We now state the main result of this section, where we denote by $2\text{-silt}(A)$ the set of silting objects in $K^{[-1,0]}(\proj\, A)$.

\begin{theorem}\cite{AIR}
Let $A$ be an algebra.
Then there is a bijection 
$$\PP(-) :  \tau\text{-tp-}A \to 2\text{-silt}(A)$$
 between the set $ \tau\text{-tp-}A$ of $\tau$-tilting pairs in $\mod\, A$ and the set $2\text{-silt}(A)$ of silting objects in $K^{[-1,0]}(\proj\, A)$, where the map is defined as
$$\PP(M,P):= \PP(M) \oplus (P \xrightarrow{} 0).$$

In particular, if $(M,P)$ is a $\tau$-rigid pair, then $\PP(M,P)$ is a presilting object in $K^{[-1,0]}(\proj\, A)$ with the same number of isomorphism classes of indecomposable direct summands.
\end{theorem}

\begin{remark}
In the Section~\ref{sec:definitions} we have discussed the the problems that one has to define the basic objects of $\tau$-tilting theory in module categories. 
However, we can see that $\tau$-tilting pairs adopt a very natural form in a category $K^{[-1,0]}(\proj\, A)$ since we don't need to distinguish between $\tau$-rigid modules and projective modules in the second entry of the pair. 
This is one of the arguments supporting the idea that $K^{[-1,0]}(\proj\, A)$ is the natural environment to study $\tau$-tilting theory. 
Even if we support this idea, we made the decision of writing this note in terms of $\tau$-tilting pairs to emphasise the relationship of $\tau$-tilting theory and classical tilting theory.
\end{remark}

\section{$g$-vectors and $\tau$-tilting theory}\label{sec:g-vectors}

In general there are many $A$-modules having the same projective presentation, which implies that $g$-vectors are in some sense ambiguous. 
However this ambiguity disappears when we restrict ourselves to $\tau$-tilting theory. 

\begin{theorem}
Let $A$ be an algebra and let $M$ and $M'$ be two $\tau$-rigid $A$-modules. 
Then $g^M=g^{M'}$ if and only if $M$ is isomorphic to $M'$.
\end{theorem}

Although the spirit of the previous result can be found already in the work of Auslander and Reiten \cite{AR1985}, the first appearance of this result stated in these terms was in the work on $2$-Calabi-Yau categories of Dehy and Keller \cite{DK}.
Later, this result was adapted to the context of $\tau$-tilting theory in the works of Adachi, Iyama and Reiten \cite{AIR} and later extended by Demonet, Iyama and Jasso in \cite{DIJ} as follows.

\begin{theorem}\cite{DIJ}
Let $A$ be an algebra and let $M$ and $M'$ be two $\tau$-rigid $A$-modules. 
Suppose that  $(g^M)_i\leq(g^{M'})_i$ for $1 \leq i \leq n$. Then $M$ is a quotient of $M'$.
In particular $g^M=g^{M'}$ if and only if $M$ is isomorphic to $M'$.
\end{theorem}

In order to state the next result we need to fix some notation. 
Given a $\tau$-tilting pair $(M,P)$ we fix a decomposition $M= \bigoplus_{i=1}^k M_i$ and $P= \bigoplus_{j= k+1}^n P_j$ of $M$ and $P$, respectively.

\begin{theorem}\label{thm:g-basis}\cite{AIR}
Let $(M,P)$ be a $\tau$-tilting pair. 
Then the set 
$$\{g^{M_1}, \dots, g^{M_k}, -g^{P_{k+1}}, \dots, -g^{P_{n}}\}$$
 form a basis of $\mathbb{Z}^n$.
\end{theorem}

\subsection{g-vectors, dimension vectors and the Euler form}

Given a finite dimensional algebra $A$, one can always associate to it a square matrix known as the Cartan matrix of the algebra as follows. 

\begin{definition}
Let $A$ be an algebra and $\{P(1), \dots, P(n)\}$ be a complete set of non-isomorphic indecomposable projective $A$-modules. 
The Cartan matrix $\mathbf{C}_A$ of $A$ is the $n \times n$ matrix 
$$\mathbf{C}_A:=([P(1)]| [P(2)]| \dots |[P(n)])$$
where the $i$-th column is the dimension vector $[P(i)]$ of $P(i)$ for $1 \leq i \leq n$. 

The Euler characteristic of $A$ is a $\mathbb{Z}$-bilinear form 
$$\langle -, - \rangle_A : K_0(A) \times K_0(A) \to \mathbb{Z}$$
defined as $\langle [M], [N] \rangle_A = [M]^{T}\mathbf{C}_A^{-1}[N]$, where $[M]$ and $[N]$ are though as column vectors.
\end{definition}

An important property of the Euler characteristic of an algebra it is that provides useful homological information, as shown in the following proposition.

\begin{proposition}
Let $A$ be an algebra of finite global dimension $s$ and let $M, N$ be two $A$-modules. 
Then 
$$\langle [M], [N] \rangle_A = \sum_{i=0}^{s} (-1)^{i}\dim_{\k}(\Ext^i_A(M,N))$$
where $\Ext^0_A(M,N)$ stands for $\Hom_A(M,N)$.
In particular, if $A$ is a hereditary algebra we have that 
$$\langle [M], [N] \rangle_A = \dim_{\k}(\Hom_A(M,N)) - \dim_{\k}(\Ext^1(M,N)).$$
\end{proposition}

\begin{remark}
The proof of the previous proposition relies onto the fact that $\{[P(1)], \dots, [P(n)]\}$ forms a basis of $\mathbb{Z}^n$ when $A$ is of finite global dimension.
\end{remark}
\smallskip

The following theorem was proven at the beginning of the 1980s by Auslander and Reiten in \cite{AR1985} but went unnoticed for several decades.
Recently, with the development of $\tau$-tilting theory this result came to light again and it is playing a key role in some of the latest developments of this theory.
To state the theorem, we denote by $\langle -, - \rangle : \mathbb{R}^n \times \mathbb{R}^n \to \mathbb{R}$ the classical dot product in $\mathbb{R}$.

\begin{theorem}\cite{AR1985}\label{thm:ARform}
Let $M$ and $N$ be modules over an algebra $A$. 
Then
$$\langle g^{M},[N]\rangle=\dim_{\k}(\Hom_A(M,N)) - \dim_{\k}(\Hom(N, \tau M)).$$
\end{theorem}

As a direct consequence of the previous result and the classical Auslander-Reiten formula we have the following corollary. 

\begin{corollary}
Let $A$ be a hereditary algebra and let $M$ and $N$ be two $A$-modules. 
Then
$$\langle g^{M},[N]\rangle = \langle [M],[N]\rangle_A.$$
\end{corollary}

Based on this last corollary, the author is of the opinion that the pairing between $g$-vectors and dimension vectors of modules is a $\tau$-tilting version of the Euler form of the algebra. 

\begin{remark}
We note that Theorem~\ref{thm:ARform} establishes a natural pairing
$$\langle -,- \rangle : K_0(K^{[-1,0]}(\proj\, A)) \times K_0(\mod\, A) \longrightarrow \mathbb{Z}$$
between $K_0(K^{[-1,0]}(\proj\, A))$ and $K_0(\mod\, A)$. 
This suggest that the extriangulated category $K^{[-1,0]}(\proj\, A)$ can be thought as a sort of "dual" category for $\mod\, A$.
We will discuss further this duality in Section~\ref{sec:wallandchamberstructure}.
\end{remark}

\subsection{$c$-vectors}

We have seen in Theorem~\ref{thm:g-basis} that the set of $g$-vectors of the indecomposable direct summands of a $\tau$-tilting pair $(M,P)$ forms a basis of $\mathbb{Z}^n$. 
This fact, together with the so-called \textit{Tropical duality} of cluster algebras \cite{FZ4, DWZ2, NZ2012}, inspired Fu to introduce in \cite{Fu2017} the notion of $c$-vector for finite dimensional algebras using $\tau$-tilting theory. 

\begin{definition}\label{def:C-mat}\cite{Fu2017}
Let $(M,P)$ be a $\tau$-tilting pair and let \sloppy $\{g^{M_{1}},\dots, g^{M_{k}},-g^{P_{k+1}},\dots, -g^{P_{n}}\}$ be the corresponding basis of $g$-vectors of $\mathbb{Z}^n$.
Define the $g$-matrix $G_{(M,P)}$ of $(M,P)$ as 
$$G_{(M,P)}=\left(g^{M_{1}} ,\dots, g^{M_{k}},-g^{P_{k+1}},\dots , -g^{P_{n}}\right)=\left( \begin{matrix}  
				(g^{M_{1}})_{1} & \dots   & (-g^{P_{n}})_{1} \\
				\vdots      & \ddots & \vdots     \\
				(g^{M_{1}})_{n}  & \dots   & (-g^{P_{n}})_{n}
		\end{matrix} \right).$$
Then the $c$-matrix $C_{(M,P)}$ of $(M,P)$ is defined as
$C_{(M,P)}=(G^{-1}_{(M,P)})^{T}.$ 
Each column of $C_{(M,P)}$ is called a \textit{$c$-vector} of $A$. 
Moreover, we say that the $i$-th column of $C_{(M,P)}$ is the $i$-th $c$-vector associated to $(M,P)$.
\end{definition}

In the same paper, Fu showed that the $c$-vectors of certain families of algebras correspond to the dimension vector of bricks. 
These results were later generalised in \cite{Treffinger2019} to every finite dimensional algebra over an algebraically closed field as follows. 

\begin{theorem}\cite{Treffinger2019}
Let $(M,P)$ be a $\tau$-tilting pair with $C$-matrix $C_{(M,P)}$. 
Then there exists a brick $B_i$ in $\mod\, A$ and $\epsilon_i \in \{0,1\}$ such that $(-1)^{\epsilon_i}[B_i] = \mathsf{c}_i$ is equal to the $i$-th $c$-vector associated to $(M,P)$.
Moreover, $[B_i]= \mathsf{c}_i$ if and only if $B_i \in \Fac\, M$. 
Dually, $-[B_i]= \mathsf{c}_i$ if and only if $B_i \in M^{\perp}$.
\end{theorem}

\begin{remark}
Note that a direct consequence of this theorem is that every $c$-vector has either only non-negative coordinates or it has only non-positive entries. 
This is usually known as \textit{sign-coherence} of $c$-vectors.
If a $c$-vectors $\mathsf{c}$ has only non-negative entries we say that $\mathsf{c}$ is a \textit{positive} $c$-vector.
Otherwise, we say that $\mathsf{c}$ is a \textit{negative} $c$-vector.
\end{remark}

In the previous result we have seen that the sign of the $c$-vectors associated to a $\tau$-tilting pair is connected with the torsion theory of $(M,P)$. 
In fact, this connection is very deep, as it can be seen in the following result. 

\begin{theorem}\cite{Treffinger2019, AsaiSemibricks}
Let $(M,P)$ be a $\tau$-tilting pair and let $\mathcal{B}_{(M,P)}$ be the set of bricks
$$\{B_i:  [B_i] = \mathsf{c}_i \text{ for some $1 \leq i \leq n$} \}.$$
Then $\mathcal{B}_{(M,P)}$ is the unique semibrick in $\mod\, A$ such that $T(\mathcal{B}_{(M,P)})= \Fac\, M$.
\end{theorem}

\begin{remark}
An interesting consequence of the precious result is that one can label the Hasse quiver $\mathbf{H}(\text{ftors-}A)$ of $\text{ftors-}A$ using (positive) $c$-vectors. 
This is compatible with the brick labelling discussed in Section~\ref{sec:MGS}.
Indeed, the labelling of $c$-vectors consists simply of taking the dimension vector of each brick in the brick labelling of $\text{ftors-}A$.
\end{remark}

In Section~\ref{sec:MGS} we spoke about the fact that we can associate to each maximal green sequence 
$$\{0\} = \T_0  \subsetneq \T_1 \subsetneq \dots \subsetneq \T_{t-1} \subsetneq \T_t = \mod\, A$$
a set $\mathcal{B}=\{B_1, \dots, B_t\}$ of bricks in $\mod\, A$. 
One can show that  the dimension vector $[B_i]$ of $B_i$ is a positive $c$-vector of $(M_i, P_i)$, where $(M_i, P_i)$ is the unique $\tau$-tilting pair such that $\T_i = \Fac\, M_i$. 
As a consequence, we can associate to every maximal green sequence a sequence $[\mathcal{B}]=\{[B_1], \dots, [B_t]\}$ of positive $c$-vectors in $\mod\, A$. 
The following result was first shown by Garver, McConville and Serhiyenko for cluster-tilted algebras in \cite{GMcS, Garver2018} and later generalised in \cite{T-HN-filt} to any algebra over an algebraically closed field.

\begin{theorem}\cite{GMcS, Garver2018, T-HN-filt}
Every maximal green sequence $\{0\} = \T_0  \subsetneq \T_1 \subsetneq \dots \subsetneq \T_{t-1} \subsetneq \T_t = \mod\, A$ is determined by its associated sequence $[\mathcal{B}]$ of $c$-vectors.
\end{theorem}


\section{The wall-and-chamber structure of an algebra}\label{sec:wallandchamberstructure}

In this final section of the paper we introduce a geometric invariant for every algebra $A$, usually known as the \textit{wall-and-chamber structure} of $A$, and we explain its relation with the $\tau$-tilting theory of the algebra.
In particular, we will see that $\tau$-tilting theory recovers much, if not all, of the stability conditions on the algebra.

\subsection{Stability conditions}

The study of stability conditions in algebraic geometry started with the introduction of geometric invariant theory by Mumford \cite{Mumford1965}.
In \cite{King1994}, King applied this theory to module categories leading to the following definition of stability conditions. 
In what follows, we always assume that the Grothendieck group $K_0(A)$ of the algebra $A$ is isomorphic to $\mathbb{Z}^n$ where the canonical basis of $K_0(A)$ is given by the set $\{[S(1)], \dots, [S(n)]\}$. 
Moreover, by abuse of notation, we identify $[M]$ with the corresponding vector of $\mathbb{Z}^n$ using the previous isomorphism.
Recall that we denote by $\langle -, - \rangle : \mathbb{R}^n \times \mathbb{R}^n \to \mathbb{R}$ the classical dot product in $\mathbb{R}$.

\begin{definition}\cite{King1994}
Let $M$ be an $A$-module with dimension vector $[M]$ and $v$ be a vector in $\mathbb{R}^n$. 
We say that $M$ is \textit{$v$-semistable} if $\langle v, [M] \rangle =0$ and $\langle v, [L]\rangle \leq 0$ for every module $L$ of $M$ different from $0$ or $M$. 
Similarly, we say that $M$ is \textit{$v$-stable} if $\langle v, [M] \rangle =0$ and $\langle v, [L]\rangle < 0 $ for every module $L$ of $M$ different from $0$ or $M$. 
\end{definition}

Later, Rudakov generalised this definition in \cite{Rudakov1997}. 
Although the Rudakov's notion of stability conditions is outwith of the scope of this note, we now state some of the algebraic consequence that he deduced in that paper. 

\begin{theorem}\cite{Rudakov1997}
Let $v$ be a vector in $\mathbb{R}^n$. 
Then the full subcategory $\mod^v_{ss} A$ of all the $v$-semistable modules in $\mod\, A$ is a wide subcategory of $\mod\, A$. 
\end{theorem}

It is not difficult ro see that the $v$-stable modules correspond to the relative simple modules in $\mod_{ss}^v A$.
The following result shows that every $v$-semistable module can be built by successive extensions of $v$-stable modules.

\begin{proposition}\cite{Rudakov1997}
Fix a vector $v$ in $\mathbb{R}^n$.
Then for every $v$-semistable module $M$ there is a filtration 
$$0 = M_0 \subset M_1 \subset \dots \subset M_{t-1} \subset M_t =M$$
where $M_i/M_{i-1}$ is a $v$-stable module. 
Moreover, all such filtrations have the same length. 
\end{proposition}

\begin{remark}
If the previous proposition reminds the reader the classical Jordan-H\"older theorem for the categories of $A$-modules, this is not a coincidence. 
We will see later that, for some specific vectors $v$, these filtrations correspond exactly with the filtrations by simple modules of an object in the module category of a smaller algebra that we can associate to $v$.
\end{remark}

In fact, one can show that stability conditions have two naturally associated torsion pairs. 
This result, first shown by Baumann, Kamnitzer and Tingley in \cite{Baumann2013}, reads as follows.

\begin{proposition}\cite{Baumann2013}
Let $v$ be a vector in $\mathbb{R}^n$. 
Then there are two torsion pairs $(\T_{\geq 0}, \F_{< 0})$ and $(\T_{>0}, \F_{\leq 0})$ where:
\begin{itemize}
	\item $\T_{\geq 0}(v):= \{X \in \mod\, A : \langle v, [Y]\rangle \geq 0 \text{ for every quotient $Y$ of $X$}\} \cup \{0\}$,
	\item $\T_{> 0}(v):= \{X \in \mod\, A : \langle v, [Y]\rangle > 0 \text{ for every quotient $Y$ of $X$}\} \cup \{0\}$,
	\item $\F_{\leq 0}(v):= \{X \in \mod\, A : \langle v, [Z]\rangle \leq 0 \text{ for every subobject $Z$ of $X$}\} \cup \{0\}$,
	\item $\F_{< 0}(v):= \{X \in \mod\, A : \langle v, [Z]\rangle < 0 \text{ for every subobject $Z$ of $X$}\} \cup \{0\}$.
\end{itemize}
Moreover the category $\mod_{ss}^v A = \T_{\geq 0}(v) \cap \F_{\leq 0}(v)$.
\end{proposition}

\subsection{The wall-and-chamber structure of an algebra}
As with many other notions in representation theory introduced since the turn of the century, one can find the origins of wall-and-chamber structures in cluster theory. 

In this case, the wall-and-chamber structure of an algebra is inspired from the notion of \textit{scattering diagrams} introduced by Gross, Hacking Keel and Kontsevich in \cite{GHKK} to study cluster algebras from a geometric perspective.
The main idea is that one can associate to any cluster algebra $\mathcal{A}$ a geometric object (its scattering diagram) that encodes much of the algebraic properties of $\mathcal{A}$.
As we have already mentioned, cluster algebras can be categorified by the representation theory of the Jacobian algebras of quivers with potentials. 
Based on this connection, Bridgeland showed in \cite{Bridgeland2016} that scattering diagrams associated to cluster algebras with an acyclic initial seed can be constructed using the stability conditions of the module category of the corresponding jacobian algebra. 
This result was later extended by Mou in \cite{Mou2019} for any cluster algebra having a green-to-red sequence. 

More presicesly, Bridgeland showed in \cite{Bridgeland2016} that one can construct a scattering diagram for any finite dimensional algebra and that this scattering diagram is isomorphic to the cluster scattering diagram of \cite{GHKK} if the algebra is hereditary.
The complete description of scattering diagrams will be skipped.
Instead, we will define now the support of the scattering diagram of an algebra, which is now called the \textit{wall-and-chamber structure} of an algebra. 
Its definition is as follows. 

\begin{definition}
The \textit{stability space} $\D(M)$ of a module $M$ is the set of vectors 
$$\D(M):= \{v \in \mathbb{R}^n : M \text{ is $v$-semistable}\}.$$
A \textit{chamber} $\Ch$ is an open connected  component of the set
$$\mathbb{R}^n \setminus \bigcup_{0 \neq M \in \mod\, A} \D(M)$$
of all vectors $v$ in $\mathbb{R}^n$ such that there is no non-zero $v$-semistable module.
A \textit{wall} is a stability space $\D(M)$ of codimension $1$, that is, a stability space $\D(M)$ such that the smallest subspace of $\mathbb{R}^n$ containing it is an hyperplane.
\end{definition}

\begin{remark}
It follows from the previous definition that the wall-and-chamber structure is completely determined by its walls (see \cite{BST2019, Asai2021}). 
Also, note that every wall is the stability space of a brick in the module category. 
However it is not true that every finite brick determines a wall.
Also, it follows from the definition that the wall-and-chamber structure of an algebra is always a fan. 
\end{remark}

\begin{remark}
We would like to emphasise that the wall-and-chamber structure of a hereditary algebra is equivalent its \textit{semi-invariant picture} introduced by Igusa, Orr, Todorov and Weymann in \cite{IOTW2009}. 
They based their construction on the notion of \textit{semi-invariant of quivers}, a notion of stability condition introduced by Schofield in \cite{Schofield1991} which is based on the Euler form of the quiver instead of of the canonical inner product of $\mathbb{R}^n$.
\end{remark}

\subsection{The $g$-vector fan of the algebra}
Recall from Section~\ref{sec:silting} that $g$-vectors correspond to the elements of the Grothendieck group $K_0(K^{[-1,0]}(\proj\, A))$ of $K^{[-1,0]}(\proj\, A)$ with respect to the basis given by the elements $\{[0 \to P(1)], \dots, [0 \to P(n)]\}$.
As such, this group is isomorphic to $\mathbb{Z}^n$. 
However, in from this section we want to relate $g$-vectors with stability conditions. 
Hence, we will see every $g$-vector $g^M$ as a vector in $\mathbb{R}^n$, where $\mathbb{R}^n = K_0(K^{[-1,0]}(\proj\, A))\otimes \mathbb{R}$ where the canonical basis is given by $\{[0 \to P(1)], \dots, [0 \to P(n)]\}$.

The aim of this subsection is to study the distribution in $\mathbb{R}^n$ of the $g$-vectors of the indecomposable $\tau$-rigid objects.
In order to do that, we start by associating a cone in $\mathbb{R}^n$ to each $\tau$-rigid pair. 
For that, recall that we are assuming that both $M$ and $P$ are basic objects that can be written as $M= \bigoplus_{i=1}^k M_i$ and $P= \bigoplus_{j=k+1}^t P_j$.

\begin{definition}
Let $A$ be an algebra and let $(M,P)$ be a $\tau$-rigid pair whose set of $g$-vectors is $$\{g^{M_1}, \dots, g^{M_k}, -g^{P_{k+1}}, \dots, -g^{P_{t}}\}.$$ 
Then we define the cone $\mathcal{C}_{(M,P)}$ to be the set
$$\mathcal{C}_{(M,P)}=\left\{\sum\limits_{i=1}\limits^{k} \alpha_ig^{M_i}-\sum\limits_{j=k+1}\limits^{t} \alpha_j g^{P_j} : \alpha_i \geq 0 \text{ for every $1\leq i\leq t$}\right\}.$$ 
Similarly, we denote by $\C^\circ_{(M,P)}$ the interior of $\C_{(M,P)}$, that is 
$$\mathcal{C}^\circ_{(M,P)}=\left\{\sum\limits_{i=1}\limits^{k} \alpha_ig^{M_i}-\sum\limits_{j=k+1}\limits^{t} \alpha_j g^{P_j} : \alpha_i > 0 \text{ for every $1\leq i\leq t$}\right\}.$$ 
\end{definition}

\begin{remark}
We know from Theorem~\ref{thm:g-basis} that $\{g^{M_1}, \dots, g^{M_k}, -g^{P_{k+1}}, \dots, -g^{P_{t}}\}$ are linearly independent. 
In particular, this implies that the cone $\C_{(M,P)}$ is of codimension $n-t$.
\end{remark}

Given two $\tau$-rigid pairs $(M_1, P_1)$ and $(M_2, P_2)$, we have by definition that $\C_{(M_1,P_1)} \cap \C_{(M_2, P_2)}$ always contains the origin $\mathbf{0} \in \mathbb{R}^n$.
Suppose that $(M,P)$ is a $\tau$-rigid pair which is a direct summand of two different $\tau$-rigid pairs $(M_1, P_1)$ and $(M_2, P_2)$. 
Then it is clear that $\C_{(M,P)} \subset \C_{(M_1, P_1)}$ and $\C_{(M,P)} \subset \C_{(M_2, P_2)}$.
In particular $\C_{(M_1,P_1)} \cap \C_{(M_2, P_2)} \neq \{\mathbf{0}\}$. 
The following result by Demonet, Iyama and Jasso in \cite{DIJ} states that this is the only way in which this can happen.

\begin{theorem}\cite{DIJ}
Let $A$ be an algebra and let $(M_1, P_1)$ and $(M_2, P_2)$ be two $\tau$-rigid pairs. 
Then $\mathcal{C}_{(M_1, P_1)} \cap \mathcal{C}_{(M_2, P_2)} \neq \{0\}$ if and only if there is a $\tau$-rigid pair $(M,P)$ which is a direct summand of both $(M_1, P_1)$ and $(M_2, P_2)$.
Moreover, if $(M,P)$ is the maximal common direct summand of $(M_1, P_1)$ and $(M_2, P_2)$  then $\mathcal{C}_{(M_1, P_1)} \cap \mathcal{C}_{(M_2, P_2)}=\mathcal{C}_{(M, P)}$.
\end{theorem}

An important consequence of the previous result is that the $g$-vectors of indecomposable $\tau$-rigid pairs have a geometrical structure known as  \textit{polyhedral fan}. 
Hence, it is common to refer the set of all $g$-vectors of indecomposable $\tau$-rigid pairs as the $g$-vector fan of the algebra. 

In particular, the $g$-vector fan of a $\tau$-tilting finite algebra has finitely many cones associated to its $\tau$-rigid pairs.
Moreover, these cones of $g$-vectors fit together very well, as was shown by Demonet, Iyama and Jasso in \cite{DIJ}.

\begin{theorem}\cite{DIJ}
Let $A$ be an algebra. 
Then $A$ is $\tau$-tilting finite if and only if 
$$\mathbb{R}^n = \bigcup_{(M,P) \text{ $\tau$-rigid}} \mathcal{C}_{(M, P)}=\left(\bigcup_{(M,P) \text{ $\tau$-rigid}} \mathcal{C}^\circ_{(M, P)}\right) \cup \{\mathbf{0}\}.$$
\end{theorem}

\subsection{$\tau$-tilting theory and stability conditions}
In Section~\ref{sec:g-vectors} we have seen that there is a natural pairing between $g$-vectors and dimension vectors that provides interesting homological information. 
For the convenience of the reader, we recall Auslander and Reiten's result here. 

\begin{theorem}\cite{AR1985}
Let $M$ and $N$ be modules over an algebra $A$. 
Then
$$\langle g^{M},[N]\rangle=\dim_{\k}(\Hom_A(M,N)) - \dim_{\k}(\Hom(N, \tau M)).$$
\end{theorem}

Using this pairing as their main tool, it was shown in \cite{BST2019} and, independently, in \cite{Yurikusa2018} that we can recover the torsion classes $\Fac\, M$ and ${}^\perp\tau M\cap P^\perp$ associated to a $\tau$-rigid pair $(M,P)$ from its $g$-vectors.
The formal statement is the following.

\begin{theorem}\cite{BST2019, Yurikusa2018}\label{thm:vss}
Let $(M,P)$ be a $\tau$-rigid pair and let $v \in \C^\circ_{(M,P)}$.
Then $\T_{>0}(v)=\Fac\, M$ and $\T_{\geq 0}(v)={}^\perp\tau M\cap P^\perp$.
In particular, $\mod_{ss}^v A = M^\perp \cap{}^\perp\tau M\cap P^\perp$.
\end{theorem}

Combining this result with Theorem~\ref{thm:tautiltred} we obtain the following corollary.

\begin{corollary}\cite{BST2019}\label{cor:vss}
Let $(M,P)$ be a $\tau$-rigid pair and let $v \in \C^\circ_{(M,P)}$.
Then the category $\mod_{ss}^v A$ of $v$-semistable modules is equivalent to the module category of the $\tau$-tilting reduction of $(M,P)$.
In particular, the number of isomorphism classes of $v$-stable modules is $n-|M|-|P|$.
\end{corollary}

\subsection{From $\tau$-tilting theory to the wall-and-chamber structure}
An important fact that follows from Corollary~\ref{cor:vss} is that the category of $v$-semistable objects is constant in a given cone $\C^\circ_{(M,P)}$.
In this subsection we explore some consequences of this result in the wall-and-chamber structure of an algebra. 

The first thing to notice is that $\mod_{ss}^v A= \{0\}$ for every vector $v\in \C^\circ_{(M,P)}$ if and only if $(M,P)$ is a $\tau$-tilting pair. 
Then, it is easy to see that $\C^\circ_{(M,P)}$ is a chamber in the wall-and-chamber structure of $A$. 
This fact was first noticed in \cite{BST2019} and it was later shown by Asai in \cite{Asai2021} that this is actually a bijection. 
The formal statement is the following. 

\begin{theorem}\cite{BST2019, Asai2021}\label{thm:chambers}
Let $(M,P)$ be a $\tau$-tilting pair. 
Then $\C^\circ_{(M,P)}$ is a chamber in the wall-and-chamber structure of $A$. 
Moreover, every chamber arises this way.
In particular, there is a one-to-one correspondence between the set of functorially finite torsion classes in $\mod\, A$ and the set of chambers of the wall-and-chamber structure of $A$.
\end{theorem}

One can also see that for every $\tau$-tilting pair $(M,P)$ the set $\C_{(M,P)} \setminus \C^\circ_{(M,P)}$ coincides with the union $\bigcup_{1\leq i \leq n} \C_{(M_i, P_i)}$ of the cones $\C_{(M_i,P_i)}$ associated to the $n$ almost-complete $\tau$-rigid pairs $\{(M_1, P_1), \dots, (M_n, P_n)\}$ that are a direct summand of $(M,P)$.
Now, it follows from Theorem~\ref{thm:vss} that there exists exactly one brick $B_i$ which is $v$-stable for every $1\leq i \leq n$ and every $v \in \C^\circ_{(M_i,P_i)}$.
Then we have built a set $\{B_1, \dots, B_n\}$ of bricks in $\mod\, A$ associated to $(M,P)$ using stability conditions. 
In the following result we show that we have encountered this set before. 

\begin{theorem}\cite{Treffinger2019}
Let $(M,P)$ be a $\tau$-tilting pair, $v\in \C^\circ_{(M,P)}$ and let $\{B_1, \dots, B_n\}$ as above. 
Then the $c$-matrix $C_{(M,P)}$ of $(M,P)$ is equal to 
$$C_{(M,P)} = (\textsf{sgn}(\langle v, [B_1]\rangle)[B_1] \quad |\quad \textsf{sgn}(\langle v, [B_2]\rangle)[B_2] \quad| \dots |\quad \textsf{sgn}(\langle v, [B_n]\rangle)[B_n]  ).$$
\end{theorem}

In particular, this implies that we can recover the semibrick associated to $(M,P)$ directly from the wall-and-chamber structure of $A$, since it is exactly the set $\{B_i : \langle v, [B_i]\rangle>0 \}$ for any $v \in \C_{(M,P)}^\circ$.

Given two chambers $\Ch_1$ and $\Ch_2$, we say that they are \textit{neighbouring} each other if $\overline{\Ch_1} \cap \overline{\Ch_2}$ is a set of codimension $1$ in $\mathbb{R}^n$.
In other words, we say that $\Ch_1$ and $\Ch_2$ are neighbours if they are separated by a wall.
Another interesting consequence of Theorem~\ref{thm:chambers} is the following.

\begin{proposition}\cite{BST2019}
Let $(M,P)$ and $(M',P')$ be two $\tau$-tilting pairs. 
Then $(M,P)$ is a mutation of $(M',P')$ if and only if the chambers $\C^\circ_{(M,P)}$ and $\C^\circ_{(M',P')}$ associated to $(M,P)$ and $(M',P')$ respectively, are neighbours.
\end{proposition}

Note that we can construct a dual graph to the wall-and-chamber structure of an algebra, where the vertices correspond to the chambers and there is an edge between two vertices if and only if the two corresponding chambers are neighbouring each other. 
Then the previous proposition is equivalent to the following. 

\begin{proposition}
The dual graph of the wall-and-chamber structure of an algebra is isomorphic to the underlying graph $\mathbf{H}(\text{ftors-}A)$, the Hasse quiver of functorially finite torsion classes in $\mod\, A$.
\end{proposition}

In particular, this result implies that crossing a wall in the wall-and-chamber of $A$ corresponds to a mutation of $\tau$-tilting pairs. 
Hence, one can realise any finite series of mutations of $\tau$-tilting pairs as paths in the wall-and-chamber structure of $A$.
Then, in particular, every maximal green sequence in $\mod\, A$ can be realised as a path in the wall-and-chamber structure of $A$, which solves the original motivating question of \cite{BST2019}.

\smallskip
Given a cone $\sigma$ in a fan $\Sigma$, the star $\mathbf{star}(\sigma)$ of $\sigma$ is the set of all the cones $\sigma'$ in $\Sigma$ having $\sigma$ as a face, that is, all the cones in $\Sigma$ such that $\sigma \subset \sigma'$.
The following result due to Asai \cite{Asai2021} shows that $\tau$-tilting reduction can be also realised at the level of wall-and-chamber structures.

\begin{theorem}\cite{Asai2021}
Let $(M,P)$ be a $\tau$-rigid pair in $\mod\, A$. 
Then the wall-and-chamber structure of the $\tau$-tilting reduction $\tilde{B}_{(M,P)}$ is isomorphic as a fan to $\mathbf{star}(\C_{(M,P)})$ of the cone $\C_{(M,P)}$ associated to $(M,P)$.
\end{theorem}

\subsection{Examples}

In this subsection we give a couple of examples illustrating the results of this section.

\begin{example}
Consider the path algebra $\mathbb{A}_2=kQ$ of the quiver $Q=\xymatrix{1\ar[r]& 2}$. Its  Auslander-Reiten quiver is as follows:

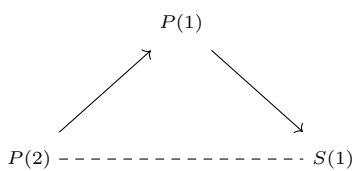
\begin{center}
			\begin{tikzpicture}[line cap=round,line join=round ,x=2.0cm,y=1.8cm]
				\clip(-1.2,-0.1) rectangle (1.2,1.1);
					\draw [->] (-0.8,0.2) -- (-0.2,0.8);
					\draw [dashed] (-0.8,0.0) -- (0.8,0.0);
					\draw [->] (0.2,0.8) -- (0.8,0.2);
				
				\begin{scriptsize}
					\draw[color=black] (-1,0) node {$P(2)$};
					\draw[color=black] (0,1) node {$P(1)$};
					\draw[color=black] (1,0) node {$S(1)$};
				\end{scriptsize}
			\end{tikzpicture}
		\end{center}

A quick calculation shows that all indecomposable $\tau$-rigid pairs in $\mod\, A$ are 
$$\left\{ (P(2), 0), (P(1), 0), (S(1), 0), (0, P(2)), (0, P(1))\right\}$$
and that the complete list of $\tau$-tilting pairs in $\mod\, A$ is
$$\left\{ (P(1)\oplus P(2), 0), (P(1)\oplus S(1), 0), (S(1), P(2)), (0, P(1)\oplus P(2)), (P(2), P(1))\right\}.$$
In Figure \ref{tauw&cA_2} we have drawn in black the wall-and-chamber structure of $A$, the cones associated to the indecomposable $\tau$-rigid pairs are in \textcolor{blue}{blue}, and the chambers associated to each $\tau$-tilting pair is indicated in \textcolor{purple}{purple}.

\begin{figure}
			\begin{center}
				\begin{tikzpicture}[line cap=round,line join=round,>=triangle 45,x=3.0cm,y=3.0cm]
					\clip(-1.5,-1.5) rectangle (1.5,1.5);
						\draw (0.,0.) -- (0.,1.5);
						\draw (0.,0.) -- (0.,-1.5);
						\draw [domain=-1.5:0.0] plot(\x,{(-0.-0.*\x)/-1.});
						\draw [domain=0.0:1.5] plot(\x,{(-0.-1.*\x)/1.});
						\draw [domain=0.0:1.5] plot(\x,{(-0.-0.*\x)/1.});
						\draw [->,line width=1.2pt,color=blue] (0.,0.) -- (1.,0.);
						\draw [->,line width=1.2pt,color=blue] (0.,0.) -- (0.,1.);
						\draw [->,line width=1.2pt,color=blue] (0.,0.) -- (-1.,0.);
						\draw [->,line width=1.2pt,color=blue] (0.,0.) -- (0.,-1.);
						\draw [->,line width=1.2pt,color=blue] (0.,0.) -- (1.,-1.);
						\draw[color=purple] (1.5,1.25) node[anchor=north east] {$\Ch_{(P(1)\oplus P(2),0)}$};
						\draw[color=purple] (-1.5,1.25) node[anchor=north west] {$\Ch_{(P(2),P(1))}$};
						\draw[color=purple] (1.55,-0.35) node[anchor=north east] {$\Ch_{(P(1)\oplus S(1), 0)}$};
						\draw[color=purple] (0.1,-1.3) node[anchor=north west] {$\Ch_{(S(1),P(2))}$};
						\draw[color=purple] (-1.5,-1.3) node[anchor=north west] {$\Ch_{(0,P(1)\oplus P(2))}$};
					\begin{scriptsize}
						\draw[color=black] (0,1.25) node[anchor= south west] {$\mathfrak{D}(S(1))$};
						\draw[color=black] (-0.2,-1.4) node {$\mathfrak{D}(S(1))$};
						\draw[color=black] (-1.25,-0.1) node {$\mathfrak{D}(S(2))$};
						\draw[color=black] (1.15,-1.4) node {$\mathfrak{D}(P(1))$};
						\draw[color=black] (1.25,0.12) node {$\mathfrak{D}(S(2))$};
						\draw[color=blue] (1.2,-0.1) node {$\C_{(P(1),0)}$};
						\draw[color=blue] (0.25,0.9) node {$\C_{(P(2),0)}$};
						\draw[color=blue] (-0.8,0.1) node {$\C_{(0,P(1))}$};
						\draw[color=blue] (0.25,-1.0) node {$\C_{(0,P(2))}$};
						\draw[color=blue] (1.25,-1.0) node {$\C_{(S(1),0)}$};
					\end{scriptsize}
				\end{tikzpicture}
			\end{center}
			\caption{Wall and chamber structure for $\mathbb{A}_2$}
            \label{tauw&cA_2}
		\end{figure}
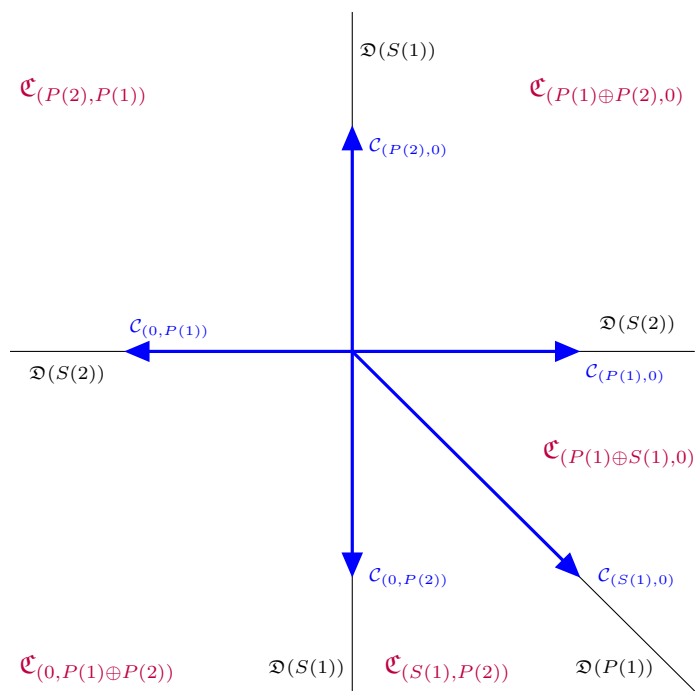
\end{example}

\begin{example}
We now take the algebra $A$ that we considered in Example~\ref{ex:supandpair}.
That is the path algebra of the quiver 
$$
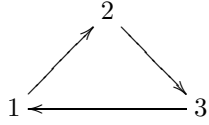
\xymatrix{
  & 2\ar[dr]& \\
  1\ar[ru] & & 3\ar[ll] }$$
modulo the second power of the ideal generated by all the arrows.
The complete list of $\tau$-tilting pairs in $\mod\, A$ can be found in Table~\ref{table:supporttautilting}.
The wall-and-chamber structure of $A$ is in $\mathbb{R}^3$ since the rank of $K_0(A)$ is equal to $3$. 
For the sake of readability, we choose to draw in Figure~\ref{fig:WC-CT} a stereographic projection of it from a vector in the first orthant. 

One can see there that all chambers in the wall-and-chamber structure of $A$ are delimited by exactly three walls. In \textcolor{red}{red} we indicate the image of the $g$-vectors of the indecomposable $\tau$-rigid pairs. 
One can obtain the $\tau$-tilting pair inducing a particular chamber from the $g$-vectors situated in the corners of that chamber. 
For instance, the chamber $\Ch_7$ is induced by the $\tau$-tilting pair $(\rep{2\\3}\oplus \rep{3}, \rep{1\\2})$.

This particular stereographic projection of the wall-and-chamber of $A$ is particularly useful to obtain the $c$-vectors associated to a given $\tau$-tilting pair, since the convexity of each of the walls indicates the sign of the corresponding $c$-vector as follows:
If the wall is convex when regarded from inside the chamber, then the $c$-vector corresponding to this wall is a negative $c$-vector and it is a positive $c$-vector otherwise. 
Using again the chamber $\Ch_7$, one can see that the $c$-vectors associated to $(\rep{2\\3}\oplus \rep{3}, \rep{1\\2})$ are $\{[0,1,0], [0,0,1], [-1, 0, -1]\}$. 
In particular, this allows us to conclude that the semibrick associated to $(\rep{2\\3}\oplus \rep{3}, \rep{1\\2})$ is $\{\rep{2} \oplus \rep{3}\}$.

Finally, note that there are $2$ walls incident to the $g$-vector $g^{\rep{3}}$, while there are three which are incident to the $g$-vector $g^{\rep{1\\2}}$.
This is due to the fact that the $\tau$-tilting reduction of $(\rep{3}, 0)$ is isomorphic to $\mathbb{A}_1 \times \mathbb{A}_1$ while that of $(\rep{1\\2}, 0)$ is isomorphic to $\mathbb{A}_2$.

\definecolor{ffzztt}{rgb}{1,0.6,0.2}
\definecolor{qqccqq}{rgb}{0,0.8,0}
\definecolor{wwqqcc}{rgb}{0.4,0,0.8}
\definecolor{yqqqyq}{rgb}{0.5019607843137255,0,0.5019607843137255}
\definecolor{qqttqq}{rgb}{0,0.2,0}
\definecolor{qqttcc}{rgb}{0,0.2,0.8}
\definecolor{ffqqtt}{rgb}{1,0,0.2}
\definecolor{ffqqqq}{rgb}{1,0,0}
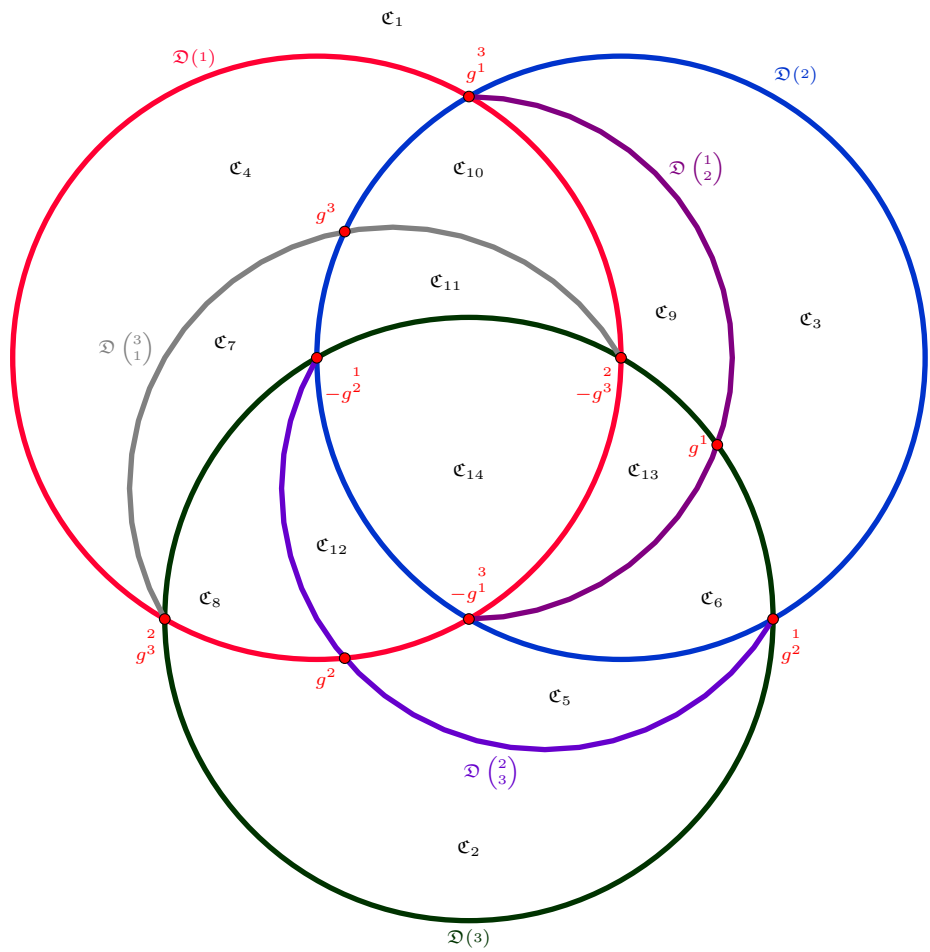
\begin{figure}
    \centering
    \begin{tikzpicture}[line cap=round,line join=round,>=triangle 45,x=1cm,y=1cm]
\clip(-7.843950349746072,-9.0) rectangle (8.613158409540794,6.839877367676253);
\draw [line width=2pt,color=ffqqtt] (-2,0) circle (4cm);
\draw [line width=2pt,color=qqttcc] (2,0) circle (4cm);
\draw [line width=2pt,color=qqttqq] (0,-3.464101615137755) circle (4cm);
\draw [shift={(0,0)},line width=2pt,color=yqqqyq]  plot[domain=-1.5707963267948966:1.5707963267948966,variable=\t]({1*3.4641016151377544*cos(\t r)+0*3.4641016151377544*sin(\t r)},{0*3.4641016151377544*cos(\t r)+1*3.4641016151377544*sin(\t r)});
\draw [shift={(1,-1.7320508075688772)},line width=2pt,color=wwqqcc]  plot[domain=2.6179938779914944:5.759586531581288,variable=\t]({1*3.4641016151377544*cos(\t r)+0*3.4641016151377544*sin(\t r)},{0*3.4641016151377544*cos(\t r)+1*3.4641016151377544*sin(\t r)});
\draw [shift={(-1,-1.7320508075688772)},line width=2pt,color=gray]  plot[domain=0.5235987755982988:3.665191429188092,variable=\t]({1*3.4641016151377544*cos(\t r)+0*3.4641016151377544*sin(\t r)},{0*3.4641016151377544*cos(\t r)+1*3.4641016151377544*sin(\t r)});
\begin{scriptsize}
\draw [fill=ffqqqq] (-2,0) circle (2.0pt) node[align=right, anchor= north west, color=red] {$-g^{\rep{1\\2}}$};
\draw [fill=ffqqqq] (2,0) circle (2.0pt) node[align=left, left, color=red, anchor=north east ] {$-g^{\rep{2\\3}}$};
\draw [fill=ffqqqq] (0,-3.464101615137755) circle (2pt) ;
\draw [fill=ffqqqq] (0,3.4641016151377544) circle (2pt) ;
\draw [fill=ffqqqq] (4,-3.4641016151377544) circle (2pt) node[align=center, above, color=red, anchor=north west] {$g^{\rep{1\\2}}$};
\draw [fill=ffqqqq] (-4,-3.4641016151377544) circle (2pt) node[align=center, above, color=red, anchor=north east] {$g^{\rep{2\\3}}$};
\draw [fill=ffqqqq] (-1.6329931618554523,1.6737265863669384) circle (2pt) node[align=center, above, color=red, anchor=south east] {$g^{\rep{3}}$};
\draw [fill=ffqqqq] (3.2659863237109037,-1.1547005383792517) circle (2pt) node[align=center, above, color=red, anchor=east] {$g^{\rep{1}}$};
\draw [fill=ffqqqq] (-1.6329931618554516,-3.9831276631254413) circle (2pt) node[align=right, above, color=red, anchor=north east] {$g^{\rep{2}}$};
\draw[color=red] (0, -3.0) node {$-g^{\rep{3\\1}}$};
\draw[color=red] (0.1, 3.9) node {$g^{\rep{3\\1}}$};
            \draw[color=ffqqtt] (-4,3.75) node[anchor=south west] {$\mathfrak{D}(\rep{1})$};
			\draw[color=qqttcc] (4.3,3.75) node {$\mathfrak{D}(\rep{2}$)};
			\draw[color=qqttqq] (0,-7.7) node {$\mathfrak{D}(\rep{3})$};
			\draw[color=yqqqyq] (3.0,2.5) node {$\mathfrak{D}\left(\rep{1\\2}\right)$};
			\draw[color=gray] (-4.5,0.12) node {$\mathfrak{D}\left(\rep{3\\1}\right)$};
			\draw[color=wwqqcc] (0.3,-5.5) node {$\mathfrak{D}\left(\rep{2\\3}\right)$};
			
			\draw (-1,4.5) node {$\Ch_{1}$};
			\draw (0,-6.5) node {$\Ch_{2}$};
			\draw (4.5,0.5) node {$\Ch_{3}$};
			\draw (-3,2.5) node {$\Ch_{4}$};
			\draw (1.2,-4.5) node {$\Ch_{5}$};
			\draw (3.2,-3.2) node {$\Ch_{6}$};
			\draw (-3.2,0.2) node {$\Ch_{7}$};
			\draw (-3.4,-3.2) node {$\Ch_{8}$};
			\draw (2.6,0.6) node {$\Ch_{9}$};
			\draw (0,2.5) node {$\Ch_{10}$};
			\draw (-0.3,1.0) node {$\Ch_{11}$};
			\draw (-1.8,-2.5) node {$\Ch_{12}$};
			\draw (2.3,-1.5) node {$\Ch_{13}$};
			\draw (0,-1.5) node {$\Ch_{14}$};
\end{scriptsize}
\end{tikzpicture}
\caption{The stereographic projection of the wall and chamber structure of $A$}
\label{fig:WC-CT}
\end{figure}
\end{example}

\subsection{Characterisations of $\tau$-tilting finite algebras}
We conclude this note by adding to Theorem~\ref{thm:tautiltingfinite} several characterisations of $\tau$-tilting finite algebras that include the notions we have discussed in this last section. 

\begin{theorem}\label{thm:tautiltingfiniteex}
Let $A$ be an algebra. 
Then the following are equivalent.
\begin{enumerate}
	\item $A$ is $\tau$-tilting finite.
	\item There are finitely many indecomposable $\tau$-rigid objects in $\mod\, A$.
	\item \cite{DIJ} There are finitely many torsion classes in $\mod\, A$.
	\item \cite{DIRRT} There are finitely many bricks in $\mod\, A$.
	\item \cite{SchrollTreffinger} The lengths of bricks in $\mod\, A$ is bounded.
	\item \cite{AsaiSemibricks} There are finitely many semibricks in $\mod\, A$.
	\item \cite{Treffinger2019} There are finitely many $c$-vectors in $\mod\, A$.
	\item \cite{MarksStovicek} There are finitely many wide subcategories in $\mod\, A$.
	\item \cite{Sentieri} Every brick in $\operatorname{Mod} A$ is a finitely presented $A$-module.
	\item \cite{DIJ} The $g$-vector fan of $A$ spans the whole $\mathbb{R}^n$.
	\item \cite{BST2019} The number of walls in the wall-and-chamber structure of $A$ is finite.
	\item \cite{BST2019} The number of chambers in the wall-and-chamber structure of $A$ is finite.
	\item \cite{Yurikusa2018} The number of different $v$-semistable subcategories of $\mod\, A$ is finite.
\end{enumerate}
\end{theorem}


\def\cprime{$'$} \def\cprime{$'$}


\vskip3mm \noindent Hipolito Treffinger\\ {\tt treffinger@imj-prg.fr}
 	
\end{document}